
\magnification=\magstep1
\vsize 23 true cm   
\hsize 15.5 true cm 
\baselineskip=13.5 pt 
\parskip=3pt 
\overfullrule=0pt

\def\BBF#1{\expandafter\edef\csname #1#1\endcsname{%
   {\mathord{\bf #1}}}}
\BBF m
\BBF p
\BBF C
\BBF P
\BBF Q
\BBF R
\BBF T
\BBF V
\BBF Z

\def\CAL#1{\expandafter\edef\csname #1\endcsname{%
   {\mathord{\cal #1}}}}
\CAL A
\CAL C
\CAL E
\CAL F
\CAL G
\CAL H
\CAL{IH}
\CAL J
\CAL K
\CAL P

\def\boule{\scriptscriptstyle\bullet}
\def\b{^{\scriptscriptstyle\bullet}}
\def\Rond#1{\mathord{\mathop{#1}\limits^{\scriptscriptstyle\circ}}}
\def\quer{\overline}
\def\mal{\mathbin{\!\cdot\!}}
\def\:{\colon}
\def\longto{\longrightarrow}
\def\into{\hookrightarrow}
\def\onto{\to\kern-.8em\to}
\def\longmapsto{\mapstochar\longto}

\def\epsilon{\varepsilon}
\def\phi{\varphi}
\def\rho{\varrho}
\def\theta{\vartheta}
\def\wideto{\;\to\;}
\def\widecong{\;\cong\;}

\def\sqr#1#2{{\vbox{\hrule height.#2pt
         \hbox{\vrule width.#2pt height#1pt \kern#1pt
               \vrule width.#2pt}
         \hrule height.#2pt}}}

\def\qed{\ifmmode\hbox{\sqr35}
     \else\hfill\sqr35\fi}

\font\bfit=cmbxti10
\font\XIIbf=cmbx12
\font\XIIrm=cmr12
\font\XIIbfM=cmbx12 scaled 1200

\def\MOPnl#1{\expandafter\edef\csname #1\endcsname{%
   {\mathop{\rm #1}\nolimits}}}
\MOPnl{cld}
\MOPnl{codim}
\MOPnl{coker}
\MOPnl{Coker}
\MOPnl{Ext}
\MOPnl{Hom}
\MOPnl{Ker}
\MOPnl{id}
\MOPnl{im}
\MOPnl{ker}
\MOPnl{lin}
\MOPnl{max}
\MOPnl{odd}
\MOPnl{or}
\MOPnl{pr}
\MOPnl{pt}
\MOPnl{rg}
\MOPnl{rk}
\MOPnl{Sp}
\MOPnl{Spec}
\MOPnl{st}
\MOPnl{supp}
\MOPnl{Tor}

\def\diagram{\def\normalbaselines{\baselineskip20pt\lineskip3pt
\lineskiplimit3pt}\matrix}
\def\mapright#1{\smash{\mathop{\hbox to 35pt{\rightarrowfill}}\limits^{#1}}}
\def\mapleft#1{\smash{\mathop{\hbox to 35pt{\leftarrowfill}}\limits^{#1}}}
\def\mapdown#1{\Big\downarrow\rlap{$\vcenter{\hbox{$\scriptstyle#1$}}$}}

\everydisplay{\abovedisplayskip 6pt plus 3 pt minus 2 pt
\belowdisplayskip\abovedisplayskip
\belowdisplayshortskip=3pt plus 3pt minus 1pt}

\font\ninerm=cmr9 
\font\ninei=cmmi9
\font\ninesy=cmsy9
\font\ninebf=cmbx9
\font\ninett=cmtt9
\font\nineit=cmti9
\font\ninesl=cmsl9

\font\sixrm=cmr6
\font\sixi=cmmi6
\font\sixsy=cmsy6
\font\sixbf=cmbx6

\def\klein{\def\rm{\fam0\ninerm}
  \textfont0=\ninerm \scriptfont0=\sixrm \scriptscriptfont0=\fiverm
  \textfont1=\ninei     \scriptfont1=\sixi \scriptscriptfont1=\fivei
  \textfont2=\ninesy     \scriptfont2=\sixsy \scriptscriptfont2=\fivesy
  \textfont\itfam=\nineit \def\it{\fam\itfam\nineit}%
  \textfont\slfam=\ninesl \def\sl{\fam\slfam\ninesl}%
  \textfont\ttfam=\ninett \def\tt{\fam\ttfam\ninett}%
  \textfont\bffam=\ninebf \scriptfont\bffam=\sixbf
   \scriptscriptfont\bffam=\fivebf \def\bf{\fam\bffam\ninebf}%
  \normalbaselineskip=10pt
  \setbox\strutbox=\hbox{\vrule height7pt depth3pt width0pt}%
  \let\sc=\sevenrm \let\big=\ninebig \normalbaselines\rm}

\newcount\annotation                                                       
 \annotation=0                                                             
  \long\def\Fussnote#1{{\baselineskip=9pt                                  
   \setbox\strutbox=\hbox{\vrule height 7pt depth 2pt width 0pt}       
    \klein\global\advance\annotation by 1                               
     \footnote{$^{\the\annotation)}$}{#1}}}

\headline{\klein\ifnum\pageno=1\else\hfill Combinatorial 
Intersection Cohomology for Fans \hfill-- \folio\ --\fi}
\nopagenumbers

\ \vskip 2 true cm

\centerline{\XIIbfM Combinatorial Intersection Cohomology for Fans}

\vskip 12 pt
\centerline{\XIIrm Gottfried Barthel, Jean-Paul Brasselet,}
\vskip 3pt
\centerline{\XIIrm Karl-Heinz Fieseler, Ludger Kaup}

\vskip 0.8 true cm
{\klein\narrower\noindent
{\bf Abstract:} Intersection cohomology $IH\b(X_\Delta;\RR)$ of a {\it 
complete\/} toric variety $X_\Delta$, associated to a fan $\Delta$ in 
$\RR^n$ and with the action of an algebraic torus $\TT \cong (\CC^*)^n$, 
is best computed using equivariant intersection cohomology 
$IH_{\TT}\b(X_\Delta)$. The reason is that $X_\Delta$ is 
$IH$-``equivariantly formal'' and equivariant intersection cohomology
provides a sheaf on $X_\Delta$, equipped with its $\TT$-invariant 
topology. An axiomatic description of that sheaf leads to the notion of 
a ``minimal  extension sheaf'' $\E\b$ on the fan $\Delta$ and a 
surprisingly simple, completely combinatorial approach which immediately 
applies to non-rational fans $\Delta$. 
These sheaves are the model for a larger class of ``pure'' sheaves, 
for which we prove a ``Decomposition Theorem''. 
For a certain class of fans 
(including fans with convex or co-convex support), called ``quasi-convex'', 
one can define a meaningful ``virtual'' intersection cohomology 
$IH\b(\Delta)$. We characterize quasi-convex fans by a purely topological 
condition on the support of their boundary fan $\partial \Delta$, and then 
deal with the question whether virtual intersection Betti numbers agree 
with the components of Stanley's generalized $h$-vector even for 
non-rational fans~$\Delta$, i.e. we try to prove that they satisfy the 
same computation algorithm. For quasi-convex fans, we prove a 
generalization of Stanley's formula realizing the intersection Poincar\'e 
polynomial of a complete toric variety in terms of local data. In order 
to show that the local data may be obtained from the virtual intersection 
cohomology of complete fans in lower dimensions, we have to assume that 
the virtual intersection cohomology of a cone $\sigma$ satisfies a certain 
vanishing condition, analoguous to the vanishing axiom for local 
intersection cohomology on the closed orbit of the affine toric variety 
$X_\sigma$ for a rational cone $\sigma$. That assumption applied to cones 
in dimension $n+1$ together with Poincar\'e duality which we show to hold 
for virtual intersection cohomology leads to a Hard Lefschetz theorem for 
polytopal fans and to the desired second step in the computation algorithm 
for virtual intersection Poincar\'e polynomials.
\par}

\vskip 0.8 true cm
{\klein
\centerline{\bf Table of Contents}
{\narrower
\medskip\noindent
Introduction                           \dotfill  2 
\par\noindent
\item{0.} Preliminaries                \dotfill  4
\item{1.} Minimal Extension Sheaves    \dotfill  8
\item{2.} Combinatorially Pure Sheaves \dotfill 13
\item{3.} Cellular Cech Cohomology     \dotfill 16
\item{4.} Poincar\'e Polynomials       \dotfill 27
\item{5.} Poincar\'e Duality           \dotfill 31
\par \vskip 3pt \noindent
References                             \dotfill 40 
\par}}

\bigskip\medskip
\centerline{\XIIbf Introduction}
\medskip\noindent
A basic combinatorial invariant of a {\it complete simplicial fan\/} 
$\Delta$ in~$\RR^n$ is its $h$-{\it vector\/}, encoding the numbers of 
cones of given dimension. By the classical {\it Dehn-Sommerville 
relations\/}, the equality $h_{i} = h_{n-i}$ holds, i.e., the vector is 
{\it palindromic\/}. 

If $\Delta$ is {\it rational\/}, then the $h$-vector admits a topological 
interpretation in terms of the associated compact $\QQ$-smooth toric 
variety~$X_{\Delta}$: By the theorem of Jurkiewicz and Danilov, the 
real\Fussnote{For ease of exposition, we use real coefficients.} 
cohomology algebra $H(X_\Delta)$ is a quotient of the Stanley-Reisner ring 
of~$\Delta$. In particular, $H(X_\Delta)$ is a combinatorial invariant 
of~$\Delta$, it ``lives'' only in even degrees, and $h_i(\Delta)$ equals 
the Betti number $b_{2i}(X_\Delta)$. Since simplicial fans are 
combinatorially equivalent to rational ones, this interpretation allows 
to apply topological results to combinatorics. Thus, the Dehn-Sommerville 
equations are just a combinatorial version of Poincar\'e duality (PD). As 
a deeper application, we mention Stanley's proof of the necessity of 
McMullen's conditions that characterize the possible $h$-vectors of 
simplicial {\it polytopal\/} fans: It involves the ``Hard'' Lefschetz 
theorem that holds since the corresponding toric variety is projective. 

In the {\it non-rational\/} case, we may ``reverse'' the theorem of 
Jurkiewicz and Danilov and take the quotient of the Stanley-Reisner ring 
as definition of a ``virtual cohomology algebra'' of the fan~$\Delta$, 
thus obtaining virtual Betti numbers $b_{2i}(\Delta)$ that coincide 
with $h_{i}(\Delta)$ for $0 \le i \le n$. \par 

In the {\it rational non-simplicial\/} case, 
using the Betti numbers of the associated 
toric variety as a definition of the $h$-vector 
no longer gives an invariant with the analoguous
properties, in fact Poincar\'e duality fails to hold,
and even worse, it is not determined by
the structure of $\Delta$ as a partially
ordered set only. Instead,
in order to get an
invariant which shares the nice
properties the classical $h$-vector has
in the simplicial case, one has to
replace singular homology with
intersection cohomology: The $i$-th
component of the {\it generalized
$h$-vector} is defined
as $h_i(\Delta):=\dim IH^{2i}(X_\Delta)$,
i.e., equals
the $2i$-th intersection Betti number
of $X_\Delta$. It satisfies Poincar\'e duality
and its components are linear functions
in the numbers of flags of cones with
prescribed sequences of dimensions. Its
computation can be done recursively using
a two step induction algorithm involving
the $g$-vector $(g_0(\sigma), \dots, g_r
(\sigma))$ of a cone $\sigma$, where
$g_i(\sigma):=\dim IH^{2i}(X_\sigma/\TT'_\sigma)$.
Here $\TT'_\sigma \subset \TT$
denotes a complementary subtorus to the stabilizer
of a point in the closed orbit of
$X_\sigma$ and $r:=[\dim \sigma / 2]-1$.
In fact, that algorithm is used to define
the generalized $h$- and the $g$-vector
also for non-rational
cones and fans, cf.~[S].
\par

In our article [Hi], we have proved
that in this situation, the r\^ole of the
Stanley-Reisner ring is played by the
$A\b:=S\b(V^*)$-module $\E\b(\Delta)$ of
global sections of a so-called ``minimal
extension sheaf" $\E\b$ on the ``fan
space" $\Delta$. (In the simplicial  case
$\E\b(\Delta)$ coincides with the
$A\b$-algebra of piecewise-polynomial
functions on $|\Delta|$, which for
complete $\Delta$ is nothing but the
Stanley-Reisner ring.) That motivates in
the non-rational case the following
definition of the {\it virtual
intersection cohomology} $IH\b(\Delta)$ of
a fan $\Delta$:  One sets
$IH\b(\Delta):=A\b/\mm \otimes_{A\b}
\E\b(\Delta)$, where $\mm:=A^{>0}$, and
hopes that for complete $\Delta$ the
components of the generalizd $h$-vector
turn out to be the virtual intersection
Betti numbers of $\Delta$, i.e.,
$h_i(\Delta) = \dim_\RR
IH^{2i}(\Delta)$.
\par

In this article we start the
investigation of the algebraic theory of
such minimal extension sheaves which
hopefully in the near future will lead to
the proof of the above interpretation of
the components of the
generalized $h$-vector. In the first
part we give the definition of minimal
extension sheaves and recall the results
of [Hi], where the virtual intersection
Betti numbers of a complete rational fan $\Delta$
are seen to equal the intersection Betti
numbers of $X_\Delta$. The second section
is devoted to combinatorially pure
sheaves over the ``fan space'' $\Delta$:
They turn out to be direct sums of simple
sheaves, which are generalized
minimal extension sheaves: To each cone
$\tau \in \Delta$ we associate a simple
pure sheaf ${}_\tau \E\b$, such that
$\E\b = {}_o \E\b$ with the zero cone
$o$, and prove a decomposition theorem
(Theorem~2.3) for pure sheaves. As a corollary, we obtain
a proof of Kalai*s conjecture for virtual intersection cohomology 
Poincar\'e polynomials, as proposed by Tom Braden, cf.\ also [BrMPh].

In the third section, we characterize ``quasi-convex'' fans, i.e., 
those fans $\Delta$ for which the $A\b$-module $\E\b (\Delta)$ is free. 
In fact, a purely $n$-dimensional fan $\Delta$ is quasi-convex 
if and only if the support of its
boundary fan $|\partial \Delta|$ is a
real homology manifold, see Theorem~3.8/9, so
in particular fans with convex or co-convex support
(i.e. $|\Delta|$ resp. $V \setminus |\Delta|$ are convex) are
quasi-convex. For rational
fans $\Delta$,
quasi-convexity is a necessary and
sufficient condition in order that the
virtual intersection cohomology agrees
with the ordinary intersection cohomology
of the associated toric variety
$X_\Delta$, i.e. $IH\b(\Delta) \cong
IH\b(X_\Delta)$. Another equivalent
reformulation of that fact is that the
intersection Betti numbers of
the associated toric variety $X_\Delta$
vanish in odd degrees. On the other hand
the freeness condition is essential in
order to have a satisfactory duality
theory both on $\E\b(\Delta)$ and
$IH\b(\Delta)=A\b/\mm \otimes_{A\b}
\E\b(\Delta)$. In fact, quasi-convexity
turns out to be equivalent to the
acyclicity of the cellular cochain
complex with coefficients in the sheaf
$\E\b$, see Theorem~3.8.
\par

The fourth
section deals with the computation of the
virtual intersection Poincar\'e
polynomials $P_\Delta:=\sum \dim
IH^{2i}(\Delta) t^{2i}$: For a
quasi-convex fan $\Delta$ the polynomial
$P_\Delta$ can be expressed in terms of
the local Poincar\'e polynomials
$P_\sigma$, see Theorem~4.3, where $P_\sigma$
denotes the virtual intersection
Poincar\'e polynomial of the fan
consisting of $\sigma$ and its proper
faces. That is a consequence of the above
mentioned acyclicity of the cellular
complex, and the fact that $IH\b(\Delta)$
and $\E\b(\Delta)$ as well as
$IH\b(\sigma)$ and $\E\b(\sigma)$ are
related by K\"unneth type formulae. In
order to get a computation algorithm for
$P_\Delta$ we have to relate
$P_\sigma$ to $P_\Lambda$, the Poincar\'e polynomial of
the ``flattened boundary fan'' $\Lambda=\Lambda_\sigma$
of $\sigma$, which in fact is a polytopal fan in $V/V_\sigma, V_\sigma
:={\rm span} (\sigma)$. For that step we need
the vanishing condition
$IH^q(\sigma)=0$ for $q \ge
\dim \sigma >0$. In the case of a
rational cone that
condition turns out to be
equivalent to the vanishing condition
for the local intersection cohomology
of $X_\sigma$ along its closed orbit; in
fact we conjecture it to hold in
general, but up to now have to state it
as a condition $V(\sigma)$, see 1.7, in the non-rational case.
The above vanishing condition together
with Poincar\'e duality (see section 5) leads to a ``Hard
Lefschetz Theorem'' for the virtual
intersection cohomology
$IH\b(\Lambda)$ of the
polytopal fan
$\Lambda_\sigma$,
see Theorem~4.6, and that theorem is behind
the description of $P_\sigma$ in terms
of $P_\Lambda$. In particular, if all
the cones in $\Delta$ satisfy the above
vanishing condition, we have
$h_i(\Delta)=\dim IH^{2i}(\Delta)$.

Finally the last section is devoted to
Poincar\'e duality: On a minimal extension sheaf $\E\b$
we define a - non-canonical - internal
intersection product $\E\b \times \E\b
\to \E\b$, that composed with an
evaluation map leads to duality
isomorphisms
$\E\b(\Delta) \cong \E\b(\Delta,
\partial \Delta)^*$ as well as
$IH\b(\Delta) \cong IH\b(\Delta,
\partial \Delta)^*$, see Theorem~5.3.

\smallskip
In order to make our results accessible to non-specialists, we 
have aimed at avoiding technical ``machinery'' and keeping the 
presentation as elementary as possible. Many essential results 
of the present article are contained in Chapters 7--10 of our 
Uppsala preprint\Fussnote{{\it ``Equivariant Intersection Cohomology 
of Toric Varieties''}, UUDM report 1998:34}; the current version 
has been announced in the note [Fi$_2$]. 
Using the formalism of derived categories, closely related work 
has been done by Tom Braden in the rational case and by Paul 
Bressler and Valery Lunts in the polytopal case. Tom Braden 
sent us a manuscript presented at the AMS meeting in 
Washington, January 2000. Even more recently, Paul Bressler and 
Valery Lunts published their ideas in the e-print [BreLu$_2$].

\smallskip
For helpful discussions, our particular
thanks go to Michel Brion,
Volker Puppe and Tom Braden.

\bigskip\medskip\goodbreak
\centerline{\XIIbf 0. Preliminaries}
\medskip\nobreak
%
\noindent
{\bf 0.A Cones and Fans:} Let $V$ be a real
vector space of dimension~$n$.  A non-zero
linear form $\alpha \: V \to \RR$ on $V$
determines the {\it upper halfspace\/}
$H_\alpha := \{v \in V; \alpha (v) \ge 0 \}$.
A (strictly convex polyhedral) {\it cone\/}
in~$V$ is a finite intersection $\sigma =
\bigcap_{i=1}^r H_{\alpha_i}$ of halfspaces
with linear forms satisfying $\bigcap_{i=1}^r
\ker\,\alpha_i =\{ 0\}$.  We let $V_\sigma :=
\sigma + (-\sigma)$ denote the linear span of
$\sigma$ in~$V$, and define $\dim \sigma :=
\dim V_\sigma$.  A cone of dimension~$d$ is
often called a {\it $d$-cone\/} for short.
\par

A cone also may be described as the set
$\sigma = \sum_{j=1}^s \RR_{\ge0}v_{j}$ of
all positive linear combinations of a finite
set of non-zero vectors~$v_{j}$ in~$V$.  In
particular, a cone spanned by a linearly
independent system of generators is called
{\it simplicial\/}.  Cones of dimension $d
\le 2$ are always simplicial; in particular,
this applies to the {\it zero cone\/} $o :=
\{0\}$ and to every {\it ray\/} (i.e., a
one-dimensional cone $\RR_{\ge0} v$).
\par

A {\it face\/} of a cone $\sigma$ is any
intersection $\tau = \sigma \cap
\ker\,\beta$, where $\beta \in V^*$ is a
linear form with $\sigma \subset H_\beta$.
We then write $\tau \preceq \sigma$ (and
$\tau \prec \sigma$ for a {\it proper\/}
face).  If in addition $\dim \tau = \dim
\sigma -1$, we call $\tau$ a {\it facet\/}
of~$\sigma$ and write $\tau \prec_1 \sigma$.
\par

\smallskip
A {\it fan\/} in $V$ is a non-empty finite
set~$\Delta$ of cones such that each face
$\tau$ of a cone $\sigma \in \Delta$ also
belongs to $\Delta$ and the intersection
$\sigma \cap \sigma'$ of two cones $\sigma,
\sigma' \in \Delta$ is a face of both,
$\sigma$ and $\sigma'$.  We say that $\Delta$
is generated by the cones $\sigma_1, \dots,
\sigma_r$, if $\Delta$ consists of all the
cones which are a face of some cone
$\sigma_i, 1 \le i \le r$.  In particular, a
given cone~$\sigma$ {\it generates\/} the fan
$\left< \sigma \right>$ consisting of
$\sigma$ and its proper faces; such a fan is
also called an {\it affine fan\/} and
occasionally is simply denoted~$\sigma$.
Moreover, we associate to~$\sigma$ its {\it
boundary fan\/} $\partial\sigma := \left<
\sigma \right> \setminus \{\sigma \}$.

Every fan is generated by the collection
$\Delta^\max$ consisting of its maximal
cones.  We define
$$
  \Delta^k := \{\sigma \in \Delta\;;\;
  \dim\sigma = k\}
  \qquad\hbox{and}\qquad \Delta^{\le k} :=
  \bigcup_{r \le k} \Delta^r \;,
$$
the latter being a subfan called the {\it
$k$-dimensional skeleton\/} (or {\it
$k$-skeleton\/} for short).  The fan $\Delta$
is called {\it purely $n$-dimensional\/} if
$\Delta^\max =\Delta^n$.  In that case, its
{\it boundary fan\/} $\partial \Delta$ is
generated by those $(n-1)$-cones that are
facets of precisely one $n$-cone in~$\Delta$.
A fan is called {\it simplicial\/} if all its
cones are simplicial; this holds if and only
if it is generated by simplicial cones.

The support $|\Delta| := \cup_{\sigma \in
\Delta} \ \sigma \subset V$ is the union of
all the cones in $\Delta$, and $\Delta$ is
called {\it complete} if and only if
$|\Delta|=V$.  We remark that the boundary
fan $\partial \Delta$ of a purely
$n$-dimensional fan~$\Delta$ is supported by
the topological boundary of~$|\Delta|$.  -- A
fan~$\Delta$ in~$V$ is called
{\it  $N$-rational\/} if there exists a
lattice (i.e., a discrete additive subgroup)
$N \subset V$ of maximal rank such that $\rho
\cap N \ne \{0\}$ for each ray $\rho \in
\Delta$.
\par

A {\it subfan\/} $\Lambda$ of a fan $\Delta$
is any subset that itself is a fan; we then
write $\Lambda \preceq \Delta$ (and $\Lambda
\prec \Delta$ if in addition~$\Lambda$ is a
{\it proper\/} subfan).  The collection of
all subfans of~$\Delta$ clearly satisfies the
axioms for the open sets of a
topology on~$\Delta$ (The empty set is admitted as a subfan).  
In the sequel, we
always endow $\Delta$ with this {\it fan
topology\/} and consider it as a topological
space, the {\it fan space}.
\par

\bigbreak
\noindent
{\bf 0.B Graded $A\b$-modules:}
In this subsection, we recall algebra results
useful for the sequel.
We denote
with $A\b$ the symmetric algebra $S\b(V^*)$
over the dual vector space $V^*$ of~$V$.  Its
elements are canonically identified with
polynomial functions on $V$.
In the case of a
rational fan $A\b$ is isomorphic to
the cohomology algebra $H\b(B\/\TT)$ of the
classifying space $B\TT$ of the complex
algebraic $n$-torus $\TT \cong (\CC^*)^n$
acting on the associated toric variety.
Motivated by that
topological considerations,  we
endow~$A\b$ with the positive even grading
determined by setting $A^{2q} := S^q(V^*)$;
in particular, $A^2 = V^*$ consists of all
linear forms on~$V$.  For a cone~$\sigma$
in~$V$, let $A\b_\sigma$ denote the algebra
$S\b(V_\sigma^*)$ with the grading as above.
The natural projection $V^* \to
V_{\sigma}^*$ extends to an epimorphism $A\b
\to A\b_\sigma$ of graded algebras.  We
usually consider the elements in $A\b_\sigma$
as functions $f \: \sigma \to \RR$; the above
epimorphism then corresponds to the
restriction of polynomial functions.
If~$\sigma$ is of dimension~$n$, then the
equality $A\b_{\sigma} = A\b$ clearly holds.
\par

For a graded $A\b$-module~$F\b$, we let
$\quer F\b$ denote its residue class module
$$
  \quer F\b \;:=\; F\b/(\mm \mal F\b)
  \widecong \RR\b \otimes_{A\b} F\b \;,
  \leqno{\rm(0.B.1)}
$$
where $\mm :=
A^{>0} \subset
A\b$ is the unique homogeneous maximal ideal
of $A\b$ and where $\RR\b := A\b/\mm = \quer
A\b$ is the field~$\RR$, considered as graded
algebra concentrated in degree zero.  Clearly
$\quer F\b$ is a graded vector space
over~$\RR$ that is finite dimensional if
$F\b$ is finitely generated over~$A\b$.  If
$F\b$ is bounded from below , then
the reverse implication holds, more precisely, a
family $(f_1, \dots, f_r)$ of homogeneous
elements in~$F\b$ generates~$F\b$ over~$A\b$
if and only if the system of residue classes
$(\quer f_1, \dots, \quer f_r)$ modulo $\mm
\mal F\b$ generates the vector
space~${\quer{F}}\b$.  In that case, we have
rk$_{A\b} F\b \le \dim \quer F\b$, with
equality holding if and only if~$F\b$ is a
free $A\b$-module.  The collection $(f_1,
\dots, f_r)$
is part of a basis of the free $A\b$-module $F\b$
over $A\b$ if and only if $(\quer f_1, \dots,
\quer f_r)$ is linearly independent
over~$\RR$.  Furthermore, every homomorphism
$\phi \: F\b \to G\b$ of graded $A\b$-modules
induces a homomorphism $\quer \phi \: {\quer
F}\b \to {\quer{G}}\b$ of graded vector
spaces that is surjective if and only
if~$\phi$ is so.  If~$F\b$ is free, then
every homomorphism $\psi \: {\quer F}\b \to
{\quer{G}}\b$ can be lifted to a homomorphism
$\phi \: F\b \to G\b$ (i.e., $\quer \phi =
\psi$ holds); if $G\b$ is free, then~$\phi$
is an isomorphism if and only if that holds
for~$\quer \phi$.\par

A finitely generated $A\b$-module $F\b$ is
{\it free\/} if and only if
$\Tor^{A\b}_1(F\b, \RR\b)=0$.  That condition
is obviously necessary, so let us show that
it is also sufficient: As we have seen above,
there is a surjection $(A\b)^m \to F\b$ where
$m := \dim \quer F\b$; let~$K\b$ be its
kernel.  Since $\Tor^{A\b}_1(F\b, \RR\b) =
0$, the exact sequence
$$
  0 \longto K\b \longto (A\b)^m
  \longto F\b \longto 0
$$
induces an exact sequence
$$
  0 \longto \quer K\b \longto (\quer A\b)^m
  \longto \quer F\b \longto 0\;.
$$
By construction, $(\quer A\b)^m \to \quer
F\b$ is an isomorphism, so we have $\quer K\b
= 0$ and thus also $K\b = 0$, i.e., $F\b
\cong (A\b)^m$ is~free.  \qed
\medskip

By means of the restriction map $A\b \to
A\b_\sigma$, an
$A\b_\sigma$-module~$F\b_\sigma$ is an
$A\b$-module, and there is a natural
isomorphism $\quer F\b_\sigma \cong
F\b_\sigma/(\mm_\sigma \mal F\b_\sigma)$.
Let us denote by $V_{\sigma}^{\perp}$ the orthogonal
complement of $V_{\sigma} \subset V$ in the
dual vector space~$V^*$,
we remark that, using the
Koszul complex for the $A\b$-module
$I(V_{\sigma}):=A\b \cdot V^\perp_\sigma
\subset A\b$, one finds a natural isomorphism
of vector spaces
$$
  \Tor^{A\b}_i(A\b_{\sigma}, \RR\b) \cong
  \Lambda^i V_{\sigma}^{\perp}
  \leqno{\rm(0.B.2)}
$$
over $\RR\b = A\b/\mm$.
 \par

\bigbreak
\noindent
{\bf 0.C Sheaves on a fan space:}
Let~$\F$ be a sheaf of
real vector
spaces on the fan space~$\Delta$.  Since the
``affine'' open sets $\langle \sigma \rangle
\preceq \Delta$ form a basis of the topology
of~$\Delta$, the sheaf~$\F$ is uniquely
determined by the collection of its values
$\F(\sigma) := \F\bigl(\left< \sigma \right>
\bigl)$ for $\sigma \in \Delta$, together
with the restriction homomorphisms
$\rho^\sigma_{\tau} \: \F(\sigma) \to
\F(\tau)$ for $\tau \preceq \sigma$.  On the
other hand, every such collection 
belongs to a sheaf~$\F$ on~$\Delta$, since an
``affine'' open subset $\langle \sigma
\rangle$ can not be covered by strictly
smaller open sets.  Furthermore, we note that
the sheaf $\F$ is flabby if and only if every
restriction map $\F (\sigma) \to \F(\partial
\sigma)$ is surjective.  -- In the same
spirit of ideas, sheaves on a fan 
occur in the work of 
Bressler and Lunts~[BreLu$_2$],
Brion~[Bri] and
McConnell~[McC].  \par

In particular, we consider the sheaf~$\A\b$
of graded algebras on~$\Delta$ given by
$\A\b(\sigma) := A\b_\sigma$, the restriction
homomorphisms $\rho^\sigma_{\tau} \:
A\b_\sigma \to A\b_\tau$ being the natural
maps $S\b(V_\sigma^*) \to S\b(V_\tau^*)$
induced by the inclusions $V_\tau \into
V_\sigma$ of a face $\tau \preceq \sigma$.
The sections $\A\b(\Lambda)$ on a subfan
$\Lambda \preceq \Delta$ constitute the
algebra of \hbox{($\Lambda$-)} piecewise
polynomial functions on $\Lambda$ in a
natural way.  \par

If $\F\b$ is a sheaf of $\A\b$-modules, then
every $\F\b(\Lambda)$ also is an
$A\b$-module, and if $\F\b(\sigma)$ is
finitely generated for every cone $\sigma \in
\Delta$, then so is $\F(\Lambda)$ for every
subfan $\Lambda \preceq \Delta$: That is an
immediate consequence of the fact that $A\b$
is a noetherian ring and of the inclusion
$\F\b(\Lambda) \subset \bigoplus_{\sigma \in
\Lambda^{\max}} \F\b(\sigma)$.  \par

For notational convenience, we often write
$$
  F\b_\Lambda := \F\b(\Lambda)
  \qquad\hbox{and}\qquad
  F\b_\sigma :=
  \F\b(\langle \sigma \rangle)\; ;
$$
more generally, for a pair of subfans
$(\Lambda, \Lambda_0)$, we define
$$
  F\b_{(\Lambda, \Lambda_0)} :=
  \ker (\rho^\Lambda_{\Lambda_{0}} \:
  F\b_\Lambda \longto F\b_{\Lambda_0})
$$
to be the submodule of sections on $\Lambda$
vanishing on $\Lambda_0$.  In particular, for
a purely $n$-dimensional fan~$\Delta$, we
obtain in that way the module
$$
  F\b_{(\Delta, \partial \Delta)} :=
  \ker (\rho^\Delta_{\partial \Delta} \:
  F\b_\Delta \longto F\b_{\partial \Delta})
$$
of sections over $\Delta$ with ``compact
supports''.\par

\smallskip\noindent
{\it Sheaves and residue class sheaves:} To a
sheaf $\F\b$ of $\A\b$-modules, we associate
the sheaf $\quer \F\b$ determined by the
assignment $\sigma \mapsto \quer F\b_\sigma =
\quer {\F\b(\sigma)}$.  This is a sheaf of
graded $\RR\b$-modules, i.e., of graded real
vector spaces. The sheaf~$\quer \F\b$ is
associated to the {\it presheaf\/} determined
by the assignment $\Lambda \mapsto
\quer{\F\b(\Lambda)}$, which in general is
not a sheaf: For a non-affine fan~$\Delta$,
the canonical homomorphism
$\quer{\F\b(\Delta)} \to \quer\F\b(\Delta)$
need not be an isomorphism.  As an example,
let~$\Delta$ be the fan describing the
projective line~$\PP_{1}$.  Then the real
vector space $\quer \A\b(\Delta) = \RR\b$ is
one-dimensional and concentrated in degree $0$, while
$\quer A\b_\Delta=
\quer{\A\b(\Delta)}
\cong H\b(\PP_{1})$ is the direct sum of two one-dimensional
weight subspaces in degree $0$ and $2$.

\bigbreak
\noindent
{\bf 0.D Fan constructions associated with a
cone:} In addition to the {\it affine
fan\/}~$\left< \sigma \right>$ and the {\it
boundary fan\/} $\partial \sigma$ associated
with a cone~$\sigma$, we need two more
constructions.  Firstly, if~$\sigma$ belongs
to a fan~$\Delta$, we consider the {\it
star\/}
$$
  \st_{\Delta}(\sigma) :=
  \{\gamma \in \Delta \;;\;
  \sigma \preceq \gamma\}
  \leqno{\rm(0.D.0)}
$$
of~$\sigma$ in~$\Delta$.  This set is not a
subfan of~$\Delta$ -- we note in passing that
it is the {\it closure\/} of the one-point
set $\{\sigma\}$ in the fan topology --, but
its image 
$$
  \Delta_\sigma \,:=\;
  p\bigl(\st_{\Delta}(\sigma)\bigr) \;=\;
  \{ p(\gamma)\;;\; \sigma \preceq \gamma \}
  \leqno{\rm(0.D.1)}
$$
under the quotient projection $p \: V \to
V/V_\sigma$ is a fan in $V/V_\sigma$, called
the {\it ``transversal fan''} of~$\sigma$
in~$\Delta$.

Secondly, if~$\sigma$ is a non-zero
cone~$\sigma$, we consider its {\it
``flattened boundary fan''}, the fan
$\Lambda_\sigma = \Lambda_\sigma (L)$ that is
obtained by projecting the boundary fan
$\partial \sigma$ onto the quotient vector
space $V_\sigma/L$, where~$L$ is a line
in~$V$ passing through the relative interior
$\buildrel {\scriptscriptstyle\circ} \over
\sigma$: If $\pi\: V_\sigma \to V_\sigma/L$
is the quotient projection, then we pose
$$
  \Lambda_\sigma := \pi (\partial \sigma) =
  \{ \pi(\tau)\,;\; \tau \prec \sigma \} \; .
  \leqno{\rm(0.D.2)}
$$
This fan is {\it complete}.  Restricting the
projection~$\pi$ to the support of $\partial
\sigma$ yields a (piecewise linear)
homeomorphism
$$
  \pi|_{\partial \sigma} \: |\partial \sigma|
  \longto |\Lambda_\sigma| = V_\sigma/L
$$
that in turn induces a homeomorphism
$\partial\sigma \to \Lambda_\sigma$
of fan spaces; in particular, the
combinatorial type of~$\Lambda_\sigma$ is
independent of the choice of~$L$.  Any linear
function $T \in A^2_\sigma$ not identically
vanishing on~$L$ provides an isomorphism $L
\buildrel \cong \over \longto \RR$;
furthermore, it gives rise to a decomposition
$V_{\sigma} = \ker(T) \oplus L$ and hence, to
an isomorphism $\ker(T) \cong V_{\sigma}/L$.
Identifying~$V_\sigma$ and $(V_\sigma/L) \times
\RR$ via these isomorphisms, we see that the
support $|\partial \sigma|$ of the boundary
fan is the graph of the strictly convex
$\Lambda_{\sigma}$-piecewise linear function
$f := T \circ (\pi|_{\partial \sigma})^{-1}
\: V_\sigma/L \to \RR$.

On the other hand, for a complete
fan~$\Lambda$ in a vector space~$W$ and a
$\Lambda$-strictly convex piecewise linear
function $f \: W \to \RR$, the convex hull of
the graph $\Gamma\!_{f}$ in $W \times \RR$ is
a cone $\gamma := \gamma^+(f)$ with
boundary~$\partial \gamma=\Gamma\!_{f}$.

\bigskip\medskip\goodbreak
\centerline{\XIIbf 1. Minimal Extension
Sheaves}
\bigskip\nobreak\noindent
The investigation of a ``virtual''
intersection cohomology theory for arbitrary
fans is {couched} in terms of a certain class
of sheaves on fans that we call {\it minimal
extension sheaves}.  In this section, we
introduce that notion and study some
elementary properties of such sheaves.
\par

\medskip\noindent
{\bf 1.1 Definition.} {\it A sheaf $\E\b$ of
graded $\A\b$-modules on the fan $\Delta$ is
called a {\it minimal extension sheaf\/}
(of~$\RR\b$) if it satisfies the following
conditions:

\smallskip{
\def\litem{\par\noindent
           \hangindent=36pt\ltextindent}
\def\ltextindent#1{\hbox
to \hangindent{#1\hss}\ignorespaces}
\litem{\bfit(N)} {\bfit Normalization:\/}
One has $E\b_o \cong A\b_o = \RR\b$ for the
zero cone~$o$.
\litem{\bfit(PF)} {\bfit Pointwise Freeness:\/}
For each cone $\sigma \in \Delta$, the
module~$E\b_{\sigma}$ is {\it free}
over~$A\b_{\sigma}$.
\litem{\bfit(LME)} {\bfit Local Minimal
Extension $\bmod\ \mm$:\/} For each cone
$\sigma \in \Delta \setminus \{ 0 \}$, the
restriction mapping
$$
  \rho_\sigma :=
\rho^{\sigma}_{\partial \sigma}
\: E\b_{\sigma} \longto
E\b _{\partial \sigma} $$
induces an isomorphism
$$
  \quer \rho_\sigma \:
  \quer{E}\b_{\sigma}
  \buildrel \cong \over \longto
  \quer{E}\b_{\partial\sigma}
$$
of graded real vector spaces.\par}
}\par

\smallskip
The above condition~(LME) implies that $\E\b$
is minimal in the set of all flabby sheaves
of graded $\A\b$-modules satisfying
conditions~(N) and~(PF), whence the name
``minimal extension sheaf'':

\medskip\noindent
{\bf 1.2 Remark.} {\it Let $\E\b$ be a
minimal extension sheaf on a fan~$\Delta$.
\item{i)} The sheaf $\E\b$ is flabby and
vanishes in odd degrees.
\item{ii)} For each subfan $\Lambda \preceq
\Delta$, the $A\b$-module $E\b_\Lambda$ is
finitely generated.  For each cone $\sigma
\in \Delta$, there is an isomorphim of graded
$A\b_{\sigma}$-modules
$$
  E\b_\sigma \widecong
  A\b_{\sigma} \otimes_{\RR} \quer E\b_\sigma \;.
  \leqno{(1.2.1)}
$$
\par}

\smallskip \noindent
{\it Proof:} (i) By the results of 0.B,
condition (LME) implies that
$\rho_\sigma$ is surjective for every cone
$\sigma \in \Delta$;  hence, 0.C asserts
flabbiness.

(ii) Let us assume that $E\b_\tau$ is
finitely generated for $\dim \tau \le k$,
then so is $E\b_\Lambda$ for every subfan
$\Lambda \preceq \Delta^{\le k}$, see 0.C.
In particular, if~$\sigma$ is a cone of
dimension $k+1$, then $E\b_{\partial \sigma}$
is finitely generated, whence $\quer
E\b_\sigma \cong \quer E\b_{\partial \sigma}$ is
finite-dimensional, and thus the free $A\b_\sigma$-module
$E\b_\sigma$ is
finitely generated. Now apply 0.C. Since $A\b$
only lives in even degrees, the obvious
$\RR\b$-splitting $F\b = F^{\rm even} \oplus
F^{\rm odd}$ of a graded $A\b$-module
actually is a decomposition into graded
$A\b$-submodules.  Hence, a finitely
generated $A\b$-module~$F\b$ vanishes in odd
degrees if and only if $\quer F\b$ does.
Thus, we may use induction on the $k$-skeletons
of $\Delta$ as above.  --- The isomorphism
(1.2.1) now is an immediate consequence of
the results quoted in (0.B) since the
$A\b_{\sigma}$-module $E\b_\sigma$ is free
and finitely generated.

\medskip

We now prove that on an arbitrary
fan~$\Delta$, a minimal extension sheaf can
be constructed recursively and that it is
unique up to isomorphism; hence, we may speak
of {\it the\/} minimal extension sheaf $\E\b
:= {}_{\Delta}\E\b$ of~$\Delta$.

\medskip
\noindent
{\bf 1.3 Proposition (Existence and
Uniqueness of Minimal Extension
Sheaves):} {\it On an arbitrary fan~$\Delta$,
there exists a minimal extension sheaf
$\E\b$; it is unique up to an isomorphism of
graded $\A\b$-modules.  More precisely, for
any two such sheaves~$\E\b$ and~$\F\b$
on~$\Delta$, every isomorphism $E\b_o \cong
F\b_o$ extends to an isomorphism $\E\b
\buildrel \cong \over \longto \F\b$ of graded
$\A\b$-modules.}\par

\smallskip
As to the uniqueness of such an extension,
see Remark~1.8, (iii).

\smallskip\noindent
{\it Proof:\/} For the {\it existence\/}, we
define the sheaf~$\E\b$ inductively on the
$k$-skeleton subfans~$\Delta^{\le k}$,
starting with $E_o\b := \RR\b$ for $k=0$.
For $k > 0$, we assume that~$\E\b$ has been
defined on~$\Delta^{<k}$; in particular,
$E\b_{\partial \sigma}$ exists for every cone
$\sigma \in \Delta^k$.  It thus suffices to
define $E\b_{\sigma}$, together with a
restriction homomorphism $E\b_{\sigma} \to
E\b_{\partial \sigma}$.  According to
(1.2.1), we set $E\b_{\sigma} \;:=\;
A\b_{\sigma} \otimes_{\RR} \quer
E\b_{\partial \sigma}$, and define the
restriction map using a $\RR\b$-linear
section $s \: \quer E\b_{\partial \sigma} \to
E\b_{\partial \sigma}$ of the residue class
map $E\b_{\partial \sigma} \to \quer
E\b_{\partial\sigma}$.  \par

For the {\it uniqueness\/} of minimal
extension sheaves up to isomorphism, we use
the same induction pattern and show how a
given isomorphism $\phi \: \E\b \to \F\b$ of
such sheaves on~$\Delta^{<k}$ may be extended
to~$\Delta^{k}$.  It suffices to verify that,
for each cone $\sigma \in \Delta^{k}$, there
is a lifting of $\phi_{\partial \sigma} \:
E\b_{\partial \sigma} \buildrel \cong \over
\longto F\b_{\partial \sigma}$ to an
isomorphism $\phi_{\sigma} \: E\b_{\sigma}
\buildrel \cong \over \longto F\b_{\sigma}$.
Using the results recalled in section 0.B,
the existence of such a lifting follows
easily from the properties of graded
$A\b_\sigma$-modules: We choose a homogeneous
basis $(e_1, \dots, e_r)$ of the free
$A\b_{\sigma}$-module~$E\b_{\sigma}$.
Since~$\F\b$ is a flabby sheaf, the images
$\phi_{\partial \sigma} (e_i|_{\partial
\sigma})$ in $F\b_{\partial \sigma}$ can be
extended to homogeneous sections $f_1, \dots,
f_r$ in $F\b_{\sigma}$ with $\deg e_j = \deg
f_j$.  The induced restriction isomorphism
$\quer F\b_{\sigma}
\buildrel\cong\over\longto \quer
F\b_{\partial\sigma}$ maps the residue
classes $\quer{f}_{1}, \dots, \quer{f}_r$ to
a basis of $\quer F\b_{\partial\sigma}$.  It
is now clear that these sections $f_1, \dots,
f_r$ form a basis of the free
$A\b_{\sigma}$-module $F\b_{\sigma}$, and
that $e_i \mapsto f_i$ defines a lifting
$\phi_{\sigma} \: E \b_{\sigma} \buildrel
\cong \over \longto F\b_{\sigma}$ of
$\phi_{\partial \sigma}$.  \qed

\medskip
Simplicial fans are easily characterized in
terms of minimal extension sheaves.

\medskip
\noindent
{\bf 1.4 Proposition:} {\it The following
conditions for a fan
$\Delta$ are equivalent:

\item{i)} $\Delta$ is simplicial,
\item{ii)} $\A\b$ is a minimal
extension sheaf on $\Delta$.}
\par

\smallskip\noindent
{\it Proof:\/} ``(ii) $\Longrightarrow$ (i)''
Assuming that~$\A\b$ is a minimal extension
sheaf, we show by induction on the
dimension~$d$ that for each cone $\sigma \in
\Delta^d$, the number $k$
of its rays equals~$d$, i.e., that~$\sigma$ is
simplicial.  This is always true for $d \le
2$.  As induction hypothesis, we assume that
the boundary fan $\partial \sigma$ is
simplicial.  On each ray of~$\sigma$, we
choose a non-zero vector $v_i$.  Then there
exist unique piecewise linear functions $f_i
\in A_{\partial \sigma}^2$ with $f_i(v_j) =
\delta_{ij}$ for $i,j = 1, \dots, k$.  As
these functions $f_1, \dots, f_k$ are clearly
linearly independent over $\RR$, we have
$\dim_{\RR} A_{\partial \sigma}^2 \ge k$.

Since $\quer A\b_\sigma = \RR\b$, we
have $\quer A^2_\sigma =0$.
Furthermore, we note that the
homogeneous  component of degree~2 in the
graded module $\mm  A_{\partial\sigma}\b$
is nothing but
$A^2 \cdot A^0_{\partial \sigma} =
A^2|_{\partial \sigma}=
A^2_{\sigma}|_{\partial\sigma}$.
Since~$\A\b$ is a minimal extension sheaf
by assumption, the induced  restriction
homomorphism $\quer A_{\sigma}\b \to
\quer A_{\partial \sigma}\b$ is an
isomorphism. We thus obtain  equalities
$$
  0 \;=\; \quer A^2_\sigma \;=\;
  \quer A_{\partial \sigma}^2
  \;=\; A_{\partial \sigma}^2 /
  (A ^2_{\sigma}|_{\partial \sigma}) \;,
$$
in particular yielding $d \le k \le \dim
A_{\partial \sigma}^2 = \dim
A^2_{\sigma}|_{\partial \sigma}$.  As
$\partial \sigma$ spans $V_\sigma$, we
further have $\dim A^2_{\sigma}|_{\partial
\sigma} = \dim A_{\sigma}^2 = d$, thus
yielding the desired result $k=d$.

``(i) $\Longrightarrow$ (ii)'': We again
proceed by induction on the dimension~$d$,
proving that for any simplicial cone~$\sigma$
with $\dim\sigma = d$ a minimal extension sheaf $\E\b$ on 
$\langle \sigma \rangle$
is naturally isomorphic to the sheaf~$\A\b$.  
The case $d=0$ being immediate,
let us first remark that a simplicial cone is
the sum $\sigma = \tau + \rho$ of
any facet $\tau \prec_{1} \sigma$ and the
remaining ray~$\rho$ that is not contained
in~$V_{\tau}$.  Using the decomposition
$V_\sigma = V_\tau \oplus V_\rho$ and the
corresponding projections $p: V_\sigma
\longrightarrow V_\tau$ and $q:V_\sigma
\longrightarrow V_\rho$ we can write
$A\b_\sigma \cong B\b_\tau \otimes_\RR
B\b_\rho$ with the subalgebras
$$
  B\b_\tau := p^*(S(V_\tau^*))
  \qquad\hbox{and}\qquad
  B\b_\rho :=q^*(S(V_\rho^*)) \;.
  \leqno{(1.4.1)}
$$
Then according to Lemma~1.5, we have
isomorphisms
$$
  E\b_{\sigma} \widecong
  A\b_{\sigma} \otimes_{B\b_{\tau}}
E\b_{\tau} \widecong
  A\b_{\sigma} \otimes_{B\b_{\tau}}
B\b_{\tau} \;=\; A\b_{\sigma}\;,
$$
as the facet~$\tau$ is simplicial and
hence we have $E\b_{\tau} \cong
B\b_{\tau}$ by induction hypothesis.\qed

\smallskip

\medskip
\noindent
{\bf 1.5 Lemma.}
{\it Assume the cone $\sigma$ is the
sum $\tau + \rho$ of a facet~$\tau$ and
a ray~$\rho$, with corresponding decompositions
$V_{\sigma} \cong V_{\tau} \oplus V_{\rho}$ and
$A\b_{\sigma} \cong B\b_\tau \otimes_\RR
B\b_\rho$ as in 1.4.1. Then the minimal extension sheaf on
$\left<\sigma\right>$ satisfies
$E\b_{\sigma} \cong A\b_{\sigma}
\otimes_{B\b_{\tau}} E\b_{\tau}$. In particular,
the restriction $E\b_\sigma \longrightarrow E\b_\tau$
induces an isomorphism $\quer E\b_{\sigma}
\cong \quer E\b_{\tau}$ of graded vector spaces.}

\smallskip\noindent
{\it Proof:} 
We use induction on $\dim \sigma$. For
$\gamma \prec \tau$ and $\hat \gamma
:=\gamma + \rho$, the induction
hypothesis yields that $E\b_{\hat \gamma} \cong
A\b_{\hat \gamma} \otimes_{B\b_\gamma}
E\b_\gamma$ with the algebra $B\b_\gamma \subset A\b_{\hat \gamma}$,
the image of $S\b(V_\gamma^*)$ in $A\b_{\hat \gamma}=S\b(V_{\hat
\gamma}^*)$ with respect to the map induced by the projection
$V_{\hat \gamma}=V_\gamma \oplus V_\rho \longrightarrow V_\gamma$.

Choosing a linear form $T \in
A^2_{ \sigma}$ vanishing on $V_\tau$, we may
write $A\b_{\sigma} = B\b_\tau [T]$,
$A\b_{\hat \gamma} = B\b_\gamma [T]$ and thus
$E\b_{\hat\gamma} \cong A\b_{\hat \gamma}
\otimes_{B\b_\gamma} E\b_\gamma
= E\b_\gamma[T]$. Then
there is an isomorphism
$E\b_{\partial \sigma} \cong
E\b_{\tau} \oplus T
E\b_{\partial \tau} [T]$ and it
suffices to check that
the restriction, which agrees with the
natural map
$$
  A\b_\sigma \otimes_{B\b_\tau} E\b_\tau
  \cong E\b_\tau [T] 
  =  E\b_{\tau} \oplus
  T E\b_{\tau} [T]
  \longto
  E\b_{\tau} \oplus
  T E\b_{\partial \tau} [T]\
$$
induces an isomorphism mod $\mm$.
It is onto, since $E\b_\tau
\to E\b_{\partial \tau}$
is. So the restriction mod $\mm$ is so
too, and it is into, since the
composition $E\b_\tau [T] \to
E\b_{\tau} \oplus T
E\b_{\partial \tau} [T] \to
E\b_\tau$ even is an isomorphism mod
$\mm$.\qed

\bigskip
If $\Delta$ is an $N$-rational fan for a
lattice $N \subset V$ of rank $n=\dim V$,
one associates to~$\Delta$ a toric
variety $X_{\Delta}$ with the action of
the
algebraic torus $\TT := N \otimes_\ZZ
\CC^* \cong  (\CC^*)^n$.
Let $IH\b_\TT(X_\Delta)$ denote the
equivariant  intersection cohomology of
$X_\Delta$ with real coefficients.
The following theorem, proved in [BBFK],
has been the starting point to
investigate minimal extension sheaves:

\medskip \noindent
{\bf 1.6 Theorem.} {\it Let $\Delta$ be a
rational fan and $\E\b$ a minimal extension
sheaf on~$\Delta$.
\item{i)} The assignment
$$
\IH\b_\TT \: \Lambda
\longmapsto IH\b_\TT (X_\Lambda)
$$
defines a sheaf on the fan space $\Delta$,
and that sheaf is a minimal extension sheaf.
\item{ii)} For each cone $\sigma \in \Delta$,
the (non-equivariant) intersection cohomology
sheaf $\IH\b$ of $X_\Delta$ is constant along
the corresponding $\TT$-orbit with stalks
isomorphic to~$\quer E\b_\sigma$.
\item{iii)} If~$\Delta$ is complete or is
affine of full dimension~$n$, then one has
$$
  IH\b(X_\Delta) \widecong \quer E\b_\Delta\;.
$$
\par}

\medskip
Statement (iii) will be generalized in
Theorem~3.8 to a considerably larger class of
rational fans that we call ``quasi-convex''.
-- For a non-zero {\it rational\/}
cone~$\sigma$, the vanishing axiom for
intersection cohomology together with
statement~(ii) yields $\quer E^q_{\sigma} =
0$ for $q \ge \dim \sigma$.  This fact turns
out to be a cornerstone in the recursive
computation of intersection Betti numbers
(see section~4).  In the non-rational case,
we have to state it as a condition; we
conjecture that it holds in general:

\smallskip

\noindent
{\bf 1.7 Vanishing Condition $\VV(\sigma)$:}
{\it A non-zero cone~$\sigma$ satisfies
the condition $\VV(\sigma)$ if
$$
  \quer E^q_{\sigma} \;=\; 0
  \quad\hbox{for}\quad q \ge \dim \sigma
  \leqno{(1.7.1)}
$$
holds. A fan~$\Delta$ satisfies the condition
$\VV(\Delta)$ if $\VV(\sigma)$ holds for each
non-zero cone $\sigma \in \Delta$.\/}

\medskip
We add some comments on that condition:  The
statements~(ii) and~(iii) in the following remark 
are not needed for
later results; in particular, the results
cited in their proof do not depend on these
statements.  -- Statement~(iii) has been
influenced by a remark of Tom Braden.

\medskip \noindent
{\bf 1.8 Remark.} i) If a fan $\Delta$ is
simplicial or rational, then condition
$\VV(\Delta)$ is satisfied.  \item{ii)}
Condition $\VV(\sigma)$ is equivalent to
$$
  E^q_{(\sigma, \partial \sigma)} \;=\;\{ 0\}
  \quad\hbox{for}\quad
  q \le \dim \sigma\;.
$$
\item{iii)} If $\Delta$ satisfies
$\VV(\Delta)$, then every homomorphism $\E\b
\to \F\b$ between minimal extension sheaves
on $\Delta$ is determined by the homomorphism
$\RR\b \cong E\b_o \to F\b_o \cong \RR\b$,
see Proposition~1.3.
\medskip

\noindent
{\it Proof:} (i) The rational case has been
mentioned above; for the simplicial case, see
Proposition~1.4.\par

(ii) We may assume $\dim\sigma = n$; hence,
the affine fan~$\left<\sigma\right>$ is
``quasi-convex'' (see Theorem~3.9).
According to Corollary~5.6, there exists an
isomorphism of abstract vector spaces $\quer
E^q_{\sigma} \cong \quer E^{2n-q}_{(\sigma,
\partial \sigma)}$. Hence $\quer E^q_{(\sigma, \partial \sigma)} 
=\{0\}$ for $q \le \dim \sigma$, whence also
$E^q_{(\sigma, \partial \sigma)} 
=\{0\}$ for $q \le \dim \sigma$, since a homogeneous base of 
$\quer E^q_{(\sigma, \partial \sigma)}$ can be lifted to a homogeneous
base of the free $A\b$-module $E^q_{(\sigma, \partial \sigma)}$. 

Now we may apply the
following remark to the finitely generated
$A\b$-module $E\b_{(\sigma, \partial
\sigma)}$: For an $A\b$-module $F\b$ which is
bounded from below one has either $F\b=0$ or,
if $r$ is minimal with $F^r \not=0$, then
$\overline F^r \cong F^r$ and $\overline
F^q=0$ for $q<r$.  Thus the claim is
immediate.  \par

(iii) We use the terminology of the proof of
Proposition~1.3:
We have to show that a homomorphism
$\phi_{\partial \sigma} \: E\b_{\partial \sigma}
\to F\b_{\partial \sigma}$
extends in a unique way to a
homomorphism $\phi_{\sigma} \:
E\b_{\sigma} \to
F\b_{\sigma}$. Statement (ii) implies that the
restrictions $E^q_\sigma
\to E^q_{\partial \sigma}$ and
$F^q_\sigma
\to F^q_{\partial \sigma}$
are isomorphisms for $q \le \dim
\sigma$. Since, as a consequence of
$\VV(\sigma)$, the $A\b$-modules
$E\b_\sigma$ and $F\b_\sigma$ can be generated
by homogeneous elements of degree $<
\dim \sigma$, the assertion follows. \qed

\bigskip\medskip\goodbreak
\centerline{\XIIbf 2. Combinatorial Pure Sheaves}
\bigskip\nobreak\noindent
In the case of a rational fan,
``the'' minimal extension
sheaf is provided by the equivariant
intersection cohomology
sheaf (see Theorem~1.6). As in
intersection cohomology, such a minimal
extension sheaf may be embedded into a class of
pure sheaves; its simple
objects are generalizations of minimal
extension sheaves. We introduce such
sheaves and prove an analogue to the
decomposition theorem in intersection cohomology.

\medskip
\noindent
{\bf 2.1 Definition:} {\it A
(combinatorially) {\bfit pure\/} sheaf on a
fan space $\Delta$ is a {\it flabby\/}
sheaf~$\F\b$ of graded $\A\b$-modules such
that, for each cone $\sigma \in \Delta$, the
$A\b_{\sigma}$-module $F\b_{\sigma}$ is {\it
finitely generated and free}.}

\medskip
\noindent
{\bf 2.2 Remark:} As a consequence of the
results in section 0.B and 0.C, we may
replace flabbiness with the following
``local'' requirement: For each cone $\sigma
\in \Delta$, the restriction homomorphism
$\rho^\sigma_{\partial\sigma} \: F\b_{\sigma}
\to F\b_{\partial \sigma}$ induces a
surjective map $\quer F\b_{\sigma} \to \quer
F\b_{\partial \sigma}$.

\smallskip
Pure sheaves are built up
from simple objects whose prototypes are
generalized minimal extension sheaves:

\medskip
\noindent
{\bf (Combinatorially) Simple Sheaves:\/ }
For each cone $\sigma \in \Delta$, we
construct inductively a ``simple''  sheaf
${}_{\sigma}\E\b$ on~$\Delta$ as
follows: For a cone  $\tau \in \Delta$
with $\dim \tau \le \dim \sigma$, we set
$$
  {}_{\sigma}E\b_{\tau} \;:=\;
  {}_{\sigma}\E\b(\tau) \;:=\;
  \cases{ A\b_{\sigma} & if $\tau = \sigma$, \cr
             0        & otherwise.\cr }
$$
Now, if ${}_{\sigma}\E\b$ has been defined
on $\Delta^{\le m}$ for some $m \ge \dim
\sigma$, then for each $\tau \in
\Delta^{m+1}$, we set
$$
  {}_{\sigma}E\b_{\tau} \;:=\;
  A\b_{\tau} \otimes_{\RR} \quer
  {}_{\sigma}\quer E\b_{\partial \tau} \;.
$$
The restriction map 
$\rho^\tau_{\partial \tau} \:
{}_{\sigma}E\b_{\tau} \to
{}_{\sigma}E\b_{\partial\tau}$ is 
induced by some homogeneous $\RR\b$-linear
section
$s:{}_{\sigma}\quer{E}\b_{\partial\tau}
\to {}_{\sigma}E\b_{\partial \tau}$
of the residue class map
${}_{\sigma}E\b_{\partial \tau} \to
{}_{\sigma}\quer E\b_{\partial\tau}$.
\smallskip

Let us collect some useful facts about
these sheaves.

\medskip\noindent
{\bf 2.2b Remark:} i) The pure sheaf $\F\b :=
{}_{\sigma}\E\b$ is determined by the
following properties of its reduction
modulo~$\mm$: \item{a)} $\quer F\b_{\sigma}
\cong \RR\b$, \item{b)} for each cone $\tau
\ne \sigma$, the reduced restriction map
$\quer F\b_{\tau} \to \quer F\b_{\partial
\tau}$ is an isomorphism.  \par\smallskip

\noindent
ii) The sheaf
${}_{\sigma}\E\b$ vanishes outside 
$\st_{\Delta}(\sigma)$ and can be obtained from
a minimal extension sheaf ${}_{\Delta_\sigma} \E\b$
on the transversal fan $\Delta_\sigma$ in the following way:
Choose a decomposition $V=V_\sigma \oplus W$, and
denote $B\b \subset A\b$ the image of $S\b((V/V_\sigma)^*)$ in $A\b$
and $B\b_\sigma$ the image of $S\b(V^*_\sigma)$ with respect to the
projection with kernel $W$. Then $A\b \cong B\b_\sigma \otimes_\RR
B\b$
and on st$(\sigma)$ we have
$$
{}_\sigma \E\b \cong 
B\b_\sigma \otimes_\RR
({}_{\Delta_\sigma} \E\b)
$$ where we identify
$\Delta_\sigma$ with st$(\sigma)$.

\noindent
iii) For the zero cone~$o$, the simple sheaf
${}_{o}\E\b$ is the minimal extension sheaf
of~$\Delta$.  \par

\noindent
iv) If $\Delta$ is a {\it rational\/} fan and
$Y \subset X_{\Delta}$ the orbit closure
associated to a cone $\sigma \in \Delta$,
then the presheaf
$$
  {}_Y\IH\b_{\TT} \;:\; \Lambda \;\mapsto\;
  IH\b_\TT(Y \cap X_{\Lambda})
$$
on~$\Delta$ actually is a sheaf isomorphic
to~${}_{\sigma}\E\b$.

\smallskip
As main result of this section, we provide a
Decomposition Formula for pure sheaves.

\medskip\noindent
{\bf 2.3 Algebraic Decomposition Theorem:}
{\it Every pure sheaf $\F\b$ on~$\Delta$
admits a direct sum decomposition
$$
  \F\b \;\cong\;\;
  \bigoplus_{\sigma \in \Delta}
  \bigl({}_{\sigma}\E\b \otimes_{\RR}
  K\b_\sigma\bigr)
$$
with $K\b_\sigma := K\b_\sigma(\F\b) :=
\ker\,(\,\quer \rho^{\sigma}_{\partial\sigma}
\: \quer F\b_{\sigma} \to
\quer F\b_{\partial\sigma})$, a finite
dimensional graded vector space.}
\smallskip

Since a finite dimensional graded vector
space $K\b$ may be uniquely written in the
form $K\b = \bigoplus
\RR\b[-\ell_{i}]^{n_{i}}$, we obtain the
``classical'' formulation
$$\F\b \cong
\bigoplus_{i}
{}_{\sigma_{i}}\E\b[-\ell_{i}]^{n_{i}}$$ of
the Decomposition Theorem.
\smallskip

\noindent
{\it Proof:\/} The following
result evidently allows an inductive
construction of such a decomposition:

{\it Given a pure sheaf~$\F\b$ and a
cone~$\sigma$ of minimal dimension with
$F\b_{\sigma} \ne 0$, there is a
decomposition $\F\b = \G\b \oplus \H\b$
as a direct sum of pure $\A\b$-submodules
$\G\b \cong _{\sigma}\E\b \otimes_{\RR}
K\b_\sigma$ and $ \H\b$, where $K\b_\sigma
= \quer F\b_\sigma$.}

Starting with $m = k:=\dim \sigma$, we construct the
decomposition recursively on each skeleton
$\Delta^{\le m}$: We set
$$
  \G\b(\tau) :=
  \cases{ F\b_\sigma \cong
               A\b_{\sigma}
                   \otimes_{\RR}
                    K\b_\sigma
               & if $\tau = \sigma$, \cr
               0 & otherwise,
              \cr } \quad\hbox{and}\quad
\H \b (\tau) := \cases{ 0 &
               if $\tau
                = \sigma $, \cr
               \F \b (\tau) & otherwise.
              \cr }
$$
We now assume that for some $m \ge k$, we
have constructed the decomposition on
$\Delta^{\le m}$. In order to extend it to
$\Delta^{\le m+1}$, it suffices again  to
extend it from the boundary fan $\partial
\tau$ of some cone  $\tau \in
\Delta^{m+1}$ to the affine
fan~$\left<\tau\right>$.
By induction hypothesis, there exists a
commutative diagram
$$
  \diagram{ & & F\b_{\tau} & \onto &
  F\b_{\partial\tau} & \cong &
  G\b_{\partial\tau} \oplus H\b_{\partial\tau} \cr
  & & \mapdown{} & & \mapdown{} & & \mapdown{} \cr
  K\b_{\tau} & \into & \quer F\b_{\tau} & \onto &
  \quer F\b_{\partial\tau} & \cong &
  \quer G\b_{\partial\tau} \oplus
\quer H\b_{\partial\tau} \cr}\ .
$$

We now choose a first decomposition
$\quer F\b_\tau = K\b_\tau
\oplus L\b \oplus M\b$ satisfying $L\b
\cong \quer G\b_{\partial \tau}$ and
$M\b \cong \quer H\b_{\partial \tau}$.
We may then lift it to a decomposition
$F\b_\tau = G\b_\tau \oplus H\b_\tau$ into free
submodules such that $\quer G\b_\tau = L\b$ and
$\quer H\b_\tau = K\b_\tau \oplus M\b$
as well as $G\b_\tau|_{\partial \tau}
= G\b_{\partial \tau}$ and
$H\b_\tau|_{\partial \tau}
= H\b_{\partial \tau}$.
\qed

\medskip\noindent
{\bf 2.4 Geometric Decomposition
Theorem:} {\it Let  $\pi : \check \Delta
\to \Delta$ be a refinement map of fans
with minimal extension sheaves
${\check \E}\b$ and $\E\b$, respectively.
Then there is  a decomposition
$$
\pi_*({\check \E}\b) \;\cong\;\; \E\b
\oplus \bigoplus_{\tau \in \Delta^{\ge
2}}  {}_{\tau}\E\b \otimes K\b_\tau $$
of $\A\b$-modules with cones
$\tau \in \Delta^{\ge 2}$ and (positively)
graded vector spaces $K\b_\tau$.}
\smallskip

\noindent
{\it Proof:\/} For an application of the
Algebraic Decomposition
Theorem~2.3, we have to verify that the
flabby sheaf $\pi_*({\check \E}\b)$  is
pure. We still need to know that the
$A\b_\sigma$-modules  $\pi_*({\check
\E}\b)(\sigma)$ are free. If~$\sigma$ is
an $n$-dimensional  cone, then the affine
fan $\left< \sigma \right>$ is
quasi-convex, see section 3.  According to
Corollary~3.11, the same holds true for
the refinement  $\check \sigma
:=\pi^{-1}( \langle \sigma \rangle )
\preceq \check \Delta$. {\it Mutatis
mutandis\/}, we may  argue along the same
lines for cones of positive codimension.
--  The fact that $\pi_*({\check \E}\b)
\cong \E\b \cong \A\b$ on  ${\Delta^{\le
1}}$ provides the condition $\dim \tau
\ge 2$, while $K^{<0}_\tau = 0$ is an
obvious consequence of the corresponding
fact  for $\pi_*({\check \E}\b)$.\qed

\bigskip
\noindent
{\bf 2.5 Corollary:} {\it Let $\pi \:
\check\Delta \to \Delta$ be a simplicial
refinement of $\Delta$.  Then the minimal
extension sheaf $\E\b$ on~$\Delta$ can be
embedded as a direct factor
into the sheaf of functions on $|\Delta|$
that are
$\check\Delta'$-piecewise polynomial.}

\smallskip
\noindent
{\it Proof:\/} According to Proposition
1.4, the sheaf~${\check\A}\b$ is a minimal
extension sheaf on~$\check \Delta$.  By
Theorem~2.4, $\E\b$ is a direct subsheaf of
$\pi_{*}({\check\A}\b)$, which is the sheaf of
functions on $|\Delta|$ that are
$\check\Delta$-piecewise polynomial.\qed

\bigskip\medskip\goodbreak
\centerline{\XIIbf 3. Cellular \v{C}ech Cohomology of Minimal 
Extension Sheaves}
\bigskip\nobreak\noindent 
In this section, our main aim is to characterize those fans~$\Delta$ 
for which the $A\b$-module $E\b_\Delta$ of global sections of a minimal 
extension sheaf $\E\b$ on~$\Delta$ is free. The principal tool is a 
``cellular'' cochain complex and the corresponding cohomology associated 
with a sheaf on a fan. The great interest in that freeness condition 
is due to the ``K\"unneth formula'' $E\b_\Delta \cong A\b
\otimes_{\RR\b} \overline E\b_\Delta$, which holds in that case. It
allows us in sections 4 and 5 to compute virtual intersection Betti
numbers
and Poincar\'e duality first on the ``equivariant'' level $E\b_\Delta$
and then to pass to ``ordinary'' (virtual) intersection cohomology
$\overline E\b_\Delta$.

The following name introduced for such fans is motivated by Theorem~3.9. 

\medskip 
\noindent 
{\bf 3.1 Definition:} {\it A fan $\Delta$ is called {\bfit quasi-convex\/} 
if the $A\b$-module $E\b_{\Delta}$ is free.}
\medskip

Obviously, quasi-convex fans are purely $n$-dimensional, i.e., each 
maximal cone in $\Delta$ is of dimension~$n$. In the rational case, 
quasi-convexity can be reformulated in terms of the associated 
toric variety:
\medskip

\noindent
{\bf 3.2 Theorem:} 
{\it A rational fan $\Delta$ is quasi-convex if and only if the 
intersection cohomology of the associated toric variety $X_\Delta$ 
vanishes in odd degrees: 
$$ 
IH^{\odd}(X_\Delta; \RR) \;:=\; 
\bigoplus_{q \ge 0} IH^{2q+1}(X_\Delta; \RR) \;=\;\{ 0 \} \; . 
$$
In that case, there exists an isomorphism $IH\b(X_\Delta) \cong 
\quer E\b_\Delta$.}
\medskip

\noindent
{\it Proof:} See Proposition 1.6 in [BBFK].\qed
\medskip

As the main tool to be used in the sequel, we now introduce the complex 
of cellular cochains on a fan with coefficients in a 
sheaf $\F$. 
\medskip

\noindent
{\bf 3.3 The cellular cochain complex.} To a fan~$\Delta$ and a 
sheaf~$\F$ of real vector spaces 
on~$\Delta$, we associate its {\it cellular cochain complex\/} 
$C\b(\Delta, \F)$: The cochain modules are  
$$
  C^k(\Delta, \F) :=  \bigoplus_{\dim \sigma = n-k} \F (\sigma) 
  \quad\hbox{for}\quad 0 \le k \le n \,. \leqno{(3.3.1)}
$$
To define the coboundary operator $\delta^k: C^k \to C^{k+1}$, we first 
fix, for each cone $\sigma \in \Delta$, an orientation $\or(\sigma)$ 
of $V_{\sigma}$ such that $\or|_{\Delta^n}$ is constant. To each facet 
$\tau \prec_{1} \sigma$, we then assign the orientation coefficient 
$\or^\sigma_{\tau}: = 1$ if the orientation of $V_\tau$, followed by 
the inward normal, coincides with the orientation of $V_\sigma$, and 
$\or^\sigma_{\tau} := -1$ otherwise. We then set
$$
  \hfil\delta (f)_{\tau} := \sum_{\sigma \succ_{1} \tau}
  \or^\sigma_{\tau} \; f_{\sigma}|_{\tau} \quad\hbox{for}\quad 
  f=(f_{\sigma}) \in C^k(\Delta, \F) \quad\hbox{and}\quad \tau \in 
  \Delta^{n-k-1}\,.\leqno{(3.3.2)} 
$$

Up to a rearrangement of indices, the complex $C\b(\Delta, \E\b)$ for
a minimal extension sheaf $\E\b$ is a {\it minimal complex\/} in the 
sense of Bernstein and Lunts. We shall come back to that at the end of 
this section. 
\par

\bigskip
More generally, we also have to consider relative cellular cochain 
complexes with respect to a subfan.

\bigskip
\noindent
{\bf 3.4 Definition:} {\it For a subfan $\Lambda$ of~$\Delta$ and a 
sheaf~$\F$ of real vector spaces 
on~$\Delta$, we set
$$
  C\b(\Delta, \Lambda; \F) := C\b(\Delta; \F)/C\b(\Lambda; \F) 
  \quad\hbox{and}\quad H^q(\Delta, \Lambda; \F) := H^q(C\b(\Delta, 
  \Lambda; \F))
$$
with the induced coboundary operator $\delta\b := \delta\b(\Delta, 
\Lambda; \F)$. If~$\Delta$ is purely $n$-dimensional and 
$\Lambda 
\preceq \partial \Delta$ a purely $n-1$-dimensional subfan
with complementary subfan $\Lambda^* \preceq \partial \Delta$ (i.e.
$\Lambda^*$ is generated by the cones in $(\partial \Delta)^{n-1}
\setminus \Lambda$), 
then the restriction of sections induces 
an augmented complex 
$$
  0 \to F_{(\Delta, \Lambda^*)}
  \buildrel\delta^{-1}\over\longto C^0(\Delta, \Lambda; \F) 
  \buildrel\delta^0\over\longto \dots  \longto C^n(\Delta, \Lambda; \F) 
  \to 0 \leqno{\tilde C\b(\Delta, \Lambda; \F):} 
$$
with cohomology $\tilde H^q(\Delta, \Lambda; \F) := 
H^q\bigl(\tilde C\b(\Delta, \Lambda; \F)\bigr)$.
}

\medskip
In fact, we need only the two cases $\Lambda= \partial \Delta$
and $\Lambda= \emptyset$, where the complementary subfan is
$\Lambda^*=\emptyset$ resp. $\Lambda^*=\partial \Delta$.
We mainly are interested in the case where~$\F$ is an $\A\b$-module. 
Then, the cohomology $\tilde H^q(\Delta, \Lambda; \F)$ is  an 
$A\b$-module. -- In the augmented situation described above, we note 
that $C^0(\Lambda; \F) = 0$ and hence $C^0(\Delta, \Lambda; \F) = 
C^0(\Delta; \F)$ holds. 

We want to compare the above cohomology in the case of the constant 
sheaf $\F = \RR$ with the usual real singular homology of a 
``spherical'' cell complex associated with a purely 
$n$-dimensional fan~$\Delta$. To that end, we fix a euclidean norm 
on~$V$ and denote with $S_V \subset V$ its unit sphere. For a subfan 
$\Lambda$ of~$\Delta$, we set 
$$
  S_\Lambda :=  |\Lambda| \cap S_V\ .
$$
For each non-zero cone~$\sigma$ in~$V$, the subset  $S_{\sigma} := 
\sigma \cap S_V$ is a  closed cell of dimension $\dim\sigma - 1$. 
Hence, the collection $(S_{\sigma})_{\sigma \in \Delta \setminus \{o\}}$ 
is a cell decomposition of $S_\Delta$, and the corresponding 
(augmented)
``homological'' complex $C_{\boule}(S_\Delta; \RR)$ of cellular chains 
with real coefficients essentially coincides with the cochain complex 
$C\b(\Delta; \RR)$: We have $C^q(\Delta; \RR) = C_{n-1-q}(S_\Delta; \RR)$ 
and $\delta^q = \partial_{n-1-q}$ for $q \le n-1$. 
 
Let us call a {\it facet-connected component} of $\Delta$ each 
purely $n$-dimensional subfan $\Delta_0$ being maximal with the 
property that every two $n$-dimensional
cones $\sigma, \sigma' \in \Lambda$ can
be joined by a chain $\sigma_0=\sigma,
\sigma_1,\dots ,\sigma_r=\sigma'$ of $n$-dimensional cones, where two
consecutive ones meet in a facet.

\medskip\noindent
{\bf 3.5 Remark:} Let $\Delta$ be a purely $n$-dimensional fan. 

\item{(i)} If $\Delta$ is complete or $n \le 1$, then 
$\tilde H\b(\Delta, \partial \Delta; \RR) = 0$.

\item{(ii)} If $\Delta$ is not complete and $n \ge 2$, then 
$$
  H^q(\Delta, \partial\Delta; \RR) \widecong 
  H_{n-1-q}(S_\Delta, S_{\partial \Delta};\RR) 
  \quad\hbox{for}\quad q > 0\;; 
$$
in particular, $H^q(\Delta, \partial\Delta; \RR) = 0$ holds for 
$q \ge n-1$.

\item{(iii)} If $s$ is the number of facet-connected components of 
$\Delta$, then $\tilde H^0(\Delta, \partial\Delta; \RR) \cong \RR^{s-1}$.
\par

\smallskip\noindent
{\it Proof:\/} The case $n \le 1$ is straightforward. For $n \ge 2$, 
the cohomology is computed via cellular homology; in the complete 
case, one has to use the fact that such a fan is facet-connected and 
that there is an isomorphism 
$$
  \hbox to \hsize{\null{(3.5.1)} \hfill 
  {$\displaystyle{ \tilde H^q(\Delta; \RR) \cong
  \tilde H_{n-1-q}(S_V;\RR)}$ 
  \quad\hbox{for $\;n \ge 2\;$ and $\;q \ge 1\;.$}
  \hfill\sqr35}}
$$

\medskip
In order to study the cellular cohomology of a flabby sheaf~$\F$ of 
real vector spaces on~$\Delta$, we want to write such a sheaf as a 
direct sum of simpler sheaves: To a cone $\sigma$ in $\Delta$, we 
associate  its {\it characteristic sheaf\/} ${}_{\sigma}\J$, i.e., 
$$
{}_{\sigma}\J (\Lambda)
 := 
\cases{ \RR & if $\sigma \in
\Lambda$
\cr \{ 0 \} & otherwise \cr}\ ,
$$
while the restriction homomorphisms are 
$id_\RR$ or $0$.

The following lemma 
is an elementary analogue of the ``Algebraic Decomposition Theorem~2.3, 
and in fact has been motivated by it.

\medskip
\noindent
{\bf 3.6 Lemma:}  {\it Every flabby sheaf~$\F$ of real vector spaces 
on~$\Delta$ admits a direct sum decomposition 
$$
  \F \widecong \bigoplus_{\sigma \in \Delta}
  {}_\sigma\J \otimes K_\sigma 
$$
with real vector spaces 
$  
  K_\sigma \widecong 
  \ker\bigl( \rho^\sigma_{\partial \sigma} \: 
  \F(\sigma) \longto \F (\partial\sigma)\bigr) \;.
$
} 

This decomposition obviously is unique up to isomorphism. 

\smallskip
\noindent
{\it Proof:\/}
The following arguments are analoguous to those in the 
proof of the Decomposition Theorem~2.3.

It clearly suffices to decompose such a flabby sheaf~$\F$ as a direct sum 
$$
  \F = \G \oplus \H
$$ 
of flabby subsheaves 
$\G$ and $\H$ such that for some cone $\sigma
\in \Delta$, we have $\G \cong {}_\sigma \J
\otimes K_\sigma$ and $\H(\sigma)=0$:  Use
induction over the number of cones $\tau
\in \Delta$, such that $\F(\tau)\ne 0$.  \par

Choose $k$ minimal
such that  there is a $k$-dimensional
cone $\sigma$ with $\F(\sigma)
\ne \{0\}$. Let
$K_\sigma := \F (\sigma)$ and define the
subsheaves $\G$ and $\H$ on the
$k$-skeleton $\Delta^{\le k}$ as follows:
$$ \G (\tau) := 
\cases{ K_\sigma &, if $\tau = \sigma$ \cr
        0 &, otherwise \cr }
$$
while
$$
\H (\tau) := 
\cases{ 0 &, if $\tau = \sigma$ \cr
        \F(\tau) &, otherwise\cr }\ .
$$
Now suppose that we already have
constructed a decomposition
$\F = \G \oplus \H$ on $\Delta^{\le m}$ 
for some $m \ge k$. Let $\tau$ be a cone of 
dimension $m+1$. In particular, we have 
a decomposition
$$
  \F (\partial \tau) = \G(\partial\tau) \oplus \H(\partial \tau)\ . 
$$
Since $\F$ is flabby, the restriction map 
$\rho^\tau_{\partial\tau} \: \F(\tau) \to \F (\partial \tau)$ is 
surjective. We can find a decomposition $\F(\tau) = U \oplus W$ 
into complementary subspaces $U, W \subset \F (\tau)$ 
such that~$\rho^\tau_{\partial\tau}$ induces an isomorphism 
$U \buildrel\cong\over\to \G (\partial\tau)$ and an 
epimorphism $W \onto \H(\partial \tau)$.
Now set $\G (\tau) := U$ and $\H(\tau) := W$. In
that manner, we can define $\G$ and $\H$
for all $(m+1)$-dimensional cones and thus
on $\Delta^{\le m+1}$. \qed

\bigskip
\noindent
Since cellular cohomology commutes with
direct sums and the tensor product with a
fixed vector space, there is an isomorphism 
$$
\tilde H\b (\Delta, \partial \Delta; \F)
\cong
\bigoplus_{\sigma \in \Delta}
\tilde H\b (\Delta, \partial \Delta;
{}_\sigma \J) \otimes K_\sigma \;,\leqno{(3.6.1)}
$$ 
so it suffices to compute
the cohomology of such a
characteristic sheaf ${}_\sigma \J$. 

\medskip 
\noindent 
{\bf 3.7 Remark:} For the characteristic sheaf ${}_\sigma \J$ of a cone 
$\sigma \in \Delta$, we have isomorphisms 
$$
  \tilde H\b(\Delta; {}_\sigma \J) \;\cong\;
  \tilde H\b(\Delta_\sigma; \RR) 
  \quad\hbox{and}\quad 
  \tilde H\b(\Delta, \partial \Delta; {}_\sigma \J) \;\cong\;
  \tilde H\b(\Delta_\sigma, \partial \Delta_\sigma; \RR) \;.
$$
with the transversal fan $\Delta_\sigma$ of $\sigma \in \Delta$.
In particular, Remark 3.5 ii) implies
$$
\tilde H^q(\Delta, \partial
\Delta; {}_\sigma \J)=0 \qquad\hbox {for}\quad 
 q \,>\, n - \dim \sigma - 2
$$
for every cone $\sigma \in \Delta$.

\medskip
We are now ready to formulate the main result of this section.

\bigskip
\noindent
{\bf 3.8 Theorem (Characterization of Quasi-Convex Fans):} {\it For a 
purely $n$-dimensional fan $\Delta$ and a minimal extension sheaf~$\E\b$ 
on it, the following statements are equivalent:
\item{i)} For each cone $\sigma \in \Delta$, we have
$$
\tilde H\b(\Delta_\sigma, \partial
\Delta_\sigma; \RR)= 0 \; .
$$
\item{ii)} We have
$$
\tilde H\b(\Delta, \partial \Delta;
\E\b)= 0 \;.
$$
\item{iii)} The fan $\Delta$ is quasi-convex, i.e., 
the $A\b$-module $E\b_\Delta:=\E\b(\Delta)$ is free.}
\par

\medskip
We put off the proof for a while, since we first want to deduce a 
{\it topological characterization\/} of quasi-convex fans. 
\medskip

\noindent
{\bf 3.9 Theorem:} {\it A  purely $n$-dimensional fan $\Delta$ is 
quasi-convex if and only if the support $|\partial \Delta|$ of its 
boundary fan is a real homology manifold. In particular, $\Delta$ is 
quasi-convex if~$\Delta$ is complete or if $S_\Delta$ is a closed 
topological $(n-1)$-cell,  e.g., if the support ~$|\Delta|$ or
the complement of the support $V \setminus |\Delta|$ are convex sets.}  
\medskip

\noindent
{\it Proof:} We note that condition (i) in Theorem 3.8 is satisfied 
for each cone $\sigma$ of dimension $n-1$ or~$n$, see 3.5 i), 
as well as for cones $\sigma \not\in \partial \Delta$, since then
$\Delta_\sigma$ is complete. 
In particular, that settles the case of a complete fan $\Delta$. In the 
non-complete case, the proof is achieved by Proposition~3.10.\qed

\medskip
To state the next result, we introduce this notation: For a 
cone~$\sigma$ in a fan~$\Delta$, we set $L_{\sigma} := 
S_{\Delta_{\sigma}} \subset (V/V_\sigma)$ and 
$\partial L_{\sigma} := S_{\partial \Delta_{\sigma}}$; in particular, 
we have $L_{o} = S_{\Delta}$. It is important to note that this 
cellular complex $L_{\sigma}$ in the $(d-1)$-sphere $S_{V/V_{\sigma}}$ 
(for $d := n - \dim\sigma$) may be identified with the {\it link\/} at 
an arbitrary point of the $(n - d - 1)$-dimensional stratum
$S_\sigma \setminus S_{\partial \sigma}$ of the stratified space $S_\Delta$.

\bigskip
\noindent
{\bf 3.10 Proposition:} {\it For a non-complete purely 
$n$-dimensional fan~$\Delta$, the following statements are 
equivalent:
\par

\item{i)} The fan $\Delta$ is quasi-convex.

\item{ii)} Each cone~$\sigma$ in $\partial \Delta$ satisfies the 
following condition:
\itemitem{(ii)$_{\sigma}$} The pair $(L_\sigma, \partial L_\sigma)$ is 
a real homology cell modulo boundary. 

\item{iii)} Each cone~$\sigma$ in $\partial \Delta$ satisfies the 
following condition:
\itemitem{(iii)$_{\sigma}$} The link~$L_\sigma$ has the real homology 
of a point.

\item{iv)} Each cone~$\sigma$ in $\partial \Delta$ satisfies the 
following condition:
\itemitem{(iv)$_{\sigma}$} The boundary of the link~$\partial L_\sigma$ 
has the real homology of a sphere of dimension $n - \dim \sigma -2$.
\par} 

\medskip
\noindent
{\it Proof:} We first show that statement~(ii) above and statement~(i) of 
Theorem~3.8 are equivalent, thus reducing the equivalence ``(i) 
$\Longleftrightarrow$ (ii)'' to Theorem~3.8. As has been remarked
above, it suffices to consider cones 
$\sigma \in (\partial \Delta)^{n-k}$ for $k \ge 2$. 
Part~(ii) of Remark~3.5 implies that 
$$
  H^q(\Delta_\sigma, \partial\Delta_\sigma; \RR) \widecong
  H_{k-1-q}(L_\sigma, \partial L_\sigma; \RR) \qquad\hbox{for $q > 0$}\; .
\leqno(3.10.0)
$$

For $q=0$, we use the equivalence
$$
H_{k-1}(L_\sigma, \partial L_\sigma;\RR) \cong
H^0(\Delta_\sigma, \partial
\Delta_\sigma; \RR) \cong \RR \iff 
\tilde H^0(\Delta_\sigma, \partial
\Delta_\sigma; \RR) =0 \;.
$$ 

In order to prove the equivalence of (ii), (iii), and~(iv), we use 
induction on~$n$. The case $n=0$ is vacuous, and in case $n=1$, it 
is trivial to check that (ii), (iii), and (iv) hold. We thus assume 
that the equivalence holds for every non-complete purely 
$d$-dimensional fan with $d \le n-1$. If we apply that to the fans
$\Delta_\sigma, \sigma \in \partial \Delta \setminus \{o\}$, we see
that the condition $ii)_\sigma$ is satisfied for every cone 
$\sigma \in \partial \Delta \setminus \{o\}$, if and only if
$iii)_\sigma$ resp. $iv)_\sigma$ is. Hence it suffices to show the
equivalence
of $ii)_o, iii)_o$ and $iv)_o$ under that assumption. We need the following

\noindent{\bfit Auxiliary Lemma:\/} {\it 
Let $L:=L_o$. If one of the conditions
$ii)_\sigma, iii)_\sigma,iv)_\sigma$ is satisfied for every 
cone $\sigma \in \partial \Delta \setminus \{o\}$, 
the inclusion of the relative 
interior $\Rond L:=L \setminus \partial L$ into $L$ induces an isomorphism 
$H_{\boule} (\Rond{L}) \cong H_{\boule} (L)$, i.e.,  
equivalently, the condition 
$$
  H_{\boule} (L,\Rond{L}) \;=\; \{0\} 
  \leqno(3.10.1)
$$ 
holds}. 

\medskip
\noindent
{\it Proof}.
For $i=-1, \dots, n-1$, we set $U_i := 
L \setminus (\partial L)_i$, where $(\partial L)_{i}$ is the 
$i$-skeleton of $\partial L = |\partial \Delta| \cap S_{V}$. By 
induction on~$i$, we show that $H_{\boule} (L, U_i)=0$ holds. 
This is evident for $i=-1$, and the case $i=n-1$ is what we have to 
prove. For the induction step, we use the homology sequence associated 
to the triple $(L,U_i,U_{i+1})$ and show 
$H_{\boule} (U_i,U_{i+1})=0$. Passing to a  
subdivision of barycentric type, we obtain an excision isomorphism between $H_{\boule} 
(U_i,U_{i+1})$ and 
$$
  H_{\boule}\Bigl(\bigcup_{\sigma \in (\partial\Delta)^{i+2}}\b 
  \st'(\Rond\sigma), \bigcup_{\sigma \in (\partial\Delta)^{i+2}}\b 
  \bigl(\st'(\Rond\sigma) \setminus \Rond\sigma \bigr)\Bigr) 
  \;=\; 
  \bigoplus_{\sigma \in (\partial\Delta)^{i+2}} 
  H_{\boule}\bigl(\st'(\Rond\sigma), 
  \st'(\Rond\sigma) \setminus \Rond\sigma \bigr)
$$
where $\st'(\Rond\sigma)$ denotes the open star of~$\Rond\sigma \cap S_{V}$ 
with respect to that subdivision of $L=S_{\Delta}$. Furthermore, there is 
a homeomorphism $\st'(\Rond\sigma) \cong \Rond{c}(L_{\sigma}) \times 
(\Rond\sigma\cap S_{V})$, where $\Rond{c}(L_{\sigma})$ denotes the open 
cone over $L_\sigma$. By the K\"unneth formula, we thus obtain an isomorphism 
$$
  H_{\boule}\bigl(\st'(\Rond\sigma), 
  \st'(\Rond\sigma) \setminus \Rond\sigma \bigr) \cong 
  H_{\boule} (\Rond{c}(L_{\sigma}),\Rond{c}(L_{\sigma})^*) \cong
  \tilde H_{\boule} (L_{\sigma})[-1]= \{0\} \;.
$$
since by the induction hypothesis $ii)_\sigma$ holds for every cone
$\sigma \in \partial \Delta \setminus \{ o \}$.
\qed
\medskip
\noindent
``$ii)_o \iff iii)_o$'' From the auxiliary lemma we obtain this 
chain of isomorphisms 
$$
H_q(L) \cong H_q(\Rond{L}) \cong 
H^{n-1-q}(S_V, S_V \setminus \Rond{L}) \cong
H^{n-1-q}(L, \partial L) \cong 
H_{n-1-q}(L, \partial L)^* \;,
$$
where the first one follows from the above lemma, the second one, from 
topological (Poincar\'e-Alexander-Lefschetz) duality, the third one is 
obtained by excision, and the fourth one is the obvious duality.
\qed

\smallskip
\noindent
``$iii)_o \Longrightarrow iv)_o$'': We 
may assume $n \ge 3$ and have to show that 
$\partial L$ has the same homology as an $(n-2)$-dimensional sphere. 
Using~(iii) together with the equivalent assumption~(ii), we have 
$\tilde H_{j-1}(\partial L) = \tilde H_j(L, \partial L) = \{0\}$ for 
$j \ne n-1$, and $H_{n-2}(\partial L) = \tilde H_{n-1}(L, \partial L) 
= \RR$. Now apply the long exact homology sequence of the pair $(L,
\partial L)$.

\smallskip
\noindent
``$iv)_o  \Longrightarrow iii)_o$'': It remains to 
verify that~$L$ has the homology of a point. We set $C := 
S^{n-1} \setminus \Rond L$ and look at the Mayer-Vietoris sequence 
$$
  \dots \to H_{q+1}(S^{n-1}) \to H_{q}(\partial L) 
  \to H_{q}(L) \oplus H_{q}(C) 
  \to H_{q}(S^{n-1}) \to H_{q-1}(\partial L) \to \dots
$$
associated to $S^{n-1} = L \cup C$. The hypothesis immediately yields 
$H_{q}(L) \oplus H_{q}(C) = 0$ for $1 \le q \le n-3$. The term 
$H_{n-1}(L) \oplus H_{n-1}(C)$ vanishes since both~$L$ and~$C$ 
are $(n-1)$-dimensional cell complexes in $S^{n-1}$ with non-empty 
boundary. The following arrow $H_{n-1}(S^{n-1}) \to H_{n-2}(\partial L)$ 
is thus an isomorphism $\RR \to \RR$. This implies that the mapping 
$H_{n-2}(L) \oplus H_{n-2}(C) \to H_{n-2}(S^{n-1})$ is injective, too, 
and that yields $H_{n-2}(L) = 0$. For $q=0$, we have a short exact 
sequence 
$$
  0 \longto \RR \longto H_{0}(L) \oplus H_{0}(C) \longto \RR \longto
  0\ ,
$$
and that yields the assertion.\qed

\medskip
As a consequence, we see that quasi-convexity of a purely 
$n$-dimensional fan depends only on the topology of its boundary: 

\bigskip \noindent    
{\bf 3.11
Corollary:} {\it Let $\Delta$ and $\Delta'$ be purely $n$-dimensional 
fans. If their boundaries have the same support $|\partial \Delta| = 
|\partial \Delta'|$, then $\Delta$ is quasi-convex if and only if 
$\Delta'$ is. 

In particular, that applies to the following special cases: 
\item{i)} $\Delta'$ is a refinement of~$\Delta$, 
\item{ii)} $\Delta$ and $\Delta'$ are ``complementary'' subfans, i.e., 
$\Delta \cup \Delta'$ is a complete fan, and $\Delta$ and $\Delta'$ 
have no $n$-dimensional cones in common.\par}

\medskip
We now come to the proof of Theorem 3.8:

\medskip \noindent 
{\bf Proof of Theorem 3.8:} For convenience, we briefly recall that 
we have to prove the equivalence of the following three statements 
for a purely $n$-dimensional fan~$\Delta$ and the minimal extension 
sheaf~$\E\b$:

{\it \item{i)}  For each cone $\sigma \in \Delta$, we have 
$\tilde H\b(\Delta_\sigma, \partial\Delta_\sigma; \RR) = 0$.
\item{ii)} We have $\tilde H\b(\Delta, \partial \Delta;\E\b) = 0$.
\item{iii)} The fan $\Delta$ is quasi-convex, i.e., 
the $A\b$-module $E\b_\Delta=\E\b(\Delta)$ is free.\par}

\noindent
``$(i) \iff (ii)$'': If we write  
$$
\E\b \cong \bigoplus_{\sigma \in \Delta}
{}_\sigma \J \otimes K_\sigma 
$$ 
according to~3.6, we obtain the following direct sum decomposition
$$
\tilde H\b(\Delta, \partial \Delta; \E\b) \cong
\bigoplus_{\sigma \in \Delta} \tilde H\b(\Delta_\sigma, \partial
\Delta_\sigma; \RR) \otimes K_\sigma
$$
according to remark 3.7 and the isomorphism 3.6.1. Hence it is
sufficient to see that none of the vector spaces
$K_\sigma \cong \ker(\rho^\sigma_{\partial \sigma}
: E\b_\sigma \longrightarrow E\b_{\partial \sigma})$ is zero:
Since $E\b_\sigma$ is a non-zero free 
$A\b_\sigma$-module and $E\b_{\partial \sigma}$ is a torsion module, 
the restriction homomorphism $\rho^\sigma_{\partial \sigma}$ never is 
injective.
\par

\noindent
``$(ii) \Longrightarrow (iii)$'': 
We shall use the abbreviations 
$$
  C^r := C^r(\Delta, \partial \Delta; \E\b)\;,\; I^r := 
{\im}\, \delta^{r-1}\;,\quad\hbox{and}\quad \Tor_k := \Tor^{A\b}_k \;.
$$ 
By downward induction on~$r$ , we show the vanishing statement 
$$
\Tor_k (I^r, \RR\b) = 0 \ {\rm for}\ k > r\;.
$$
In particular, the $A\b$-module $I^0 = E\b_\Delta$
satisfies $\Tor_1(I^0,\RR\b)=0$ and thus is free according to~(0.B). 
\par

Obviously the above statement holds for $r=n+1$. Since $C\b$ is acyclic 
by assumption~(ii), the sequences 
$$
0 \longto I^r \longto C^r \longto I^{r+1}
\longto 0
$$
are exact and thus induce exact sequences
$$
  \Tor_{k+1} (I^{r+1}, \RR\b) 
  \longto \Tor_k (I^r,\RR\b)
  \longto \Tor_k (C^r, \RR\b) \; .
$$
The last term vanishes for $k > r$: The module $C^r =
\bigoplus_{\dim\sigma = n-r} E\b_\sigma$ actually is a direct sum 
of shifted modules $A\b_\sigma$, hence $\Tor_k(C^r, \RR\b)=0$ for 
$k>r$, cf.\ 0.B.1. Since by induction hypothesis, also the first term 
vanishes, so does the second one.
\par

\noindent
``$(iii) \Longrightarrow (ii)$''  
In addition to the above, we use the abbreviations 
$$ 
  K^r := {\ker}\, \delta^r \quad\hbox{and}\quad 
  \tilde H^r := \tilde H^r(\Delta, \partial \Delta; \E\b) 
  = K^r/I^r \;. 
$$
We have to prove that $\tilde H^r = 0$ holds for 
all~$r$. We choose 
an increasing sequence of subspaces
$V_0 := 0 \subset V_1 \subset \dots \subset
V_n := V$ such that $V = V_r \oplus
V_\sigma$ holds for all $\sigma \in
\Delta^{n-r}$. Then the algebras $B\b_r :=
S\b((V/V_r)^*)$ form a decreasing sequence 
of subalgebras of~$A\b$, and for all 
cones $\sigma \in \Delta^{n-r}$, there is 
an isomorphism $B\b_r \cong A\b_{\sigma} $ 
induced from the composed mapping 
$V_\sigma \to V \to V/V_r$. In particular, 
each~$C^r = \bigoplus_{\sigma \in 
\Delta^{n-r}} E\b_{\sigma}$ is a free 
$B\b_r$-module. 

We now choose linear forms $T_1, \dots, T_n 
\in A^2$ such that 
$B\b_{r} = \RR[T_1, \dots, T_{n-r}]$. 
By induction on $r$, we shall prove:
$$
  \tilde H^q \;=\; 0 \quad\hbox{for}\quad q < r \;,
  \quad\hbox{and $\; I^r\;$ is a free $B\b_r$-module.}
$$
Since 
$I^0 = E\b_\Delta$, the assertion is evident 
for $r=0$. So let us proceed from~$r$ to 
$r+1$. The vanishing of~$\tilde H^r$ is a 
consequence of the fact that its support in 
$\Spec(B\b_r)$ is too small: According to 
Lemma~3.13 below, the support of $\tilde H^r$ 
in $\Spec(A\b)$ has codimension at least 
$r+2$. Thus, as $B\b_r$-module, its support 
in $\Spec(B\b_r)$ has codimension at 
least~$2$. Using the exact sequence
$$
0 \longto I^r \longto K^r \longto \tilde
H^r \longto 0 \;,
$$
the vanishing $\tilde H^r = 0$ then follows 
from Lemma~3.12.

\smallskip
It remains to prove that $I := I^{r+1}$ is a
free module over $B\b := B\b_{r+1}$. By 0.B, 
this is equivalent to 
$$
  \Tor_1^{B\b}(I, \RR) = 0 \;.\leqno{(3.8.1)}
$$ 
Recall that $B\b_r = B\b[T]$ where $T := T_{n-r}$. 
Thus, the formula
$$ 
\Tor_k^{B\b}(I, \RR) \cong
\Tor_k^{B\b[T]}(I, \RR[T]) \leqno{(3.8.2)}
$$
provides the bridge to the induction hypothesis on 
the previous level~$r$. Now the exact sequence
$$
0 \longto \RR[T] 
\buildrel \mu \over
\longto \RR[T]
\longto \RR 
\longto 0\; .\leqno{(3.8.3)}
$$
where $\mu$ is multiplication with $T$, yields an exact sequence 
$$
\Tor_2^{B\b[T]}(I, \RR) \longto
\Tor_1^{B\b[T]}(I, \RR[T])
\buildrel \theta \over \longto
\Tor_1^{B\b[T]}(I, \RR[T])\; .\leqno{(3.8.4)}
$$
The homomorphism~$\theta$ is {\it injective\/}: The vector space 
$\Tor_2^{B\b[T]}(I, \RR) \cong \Tor_2^{B\b_r}(I, \RR)$ vanishes 
since we already know that $\tilde H^r=\{0\}$ and thus the exact sequence
$$
  0 \longto I^r \longto C^r \longto I
  \longto 0 \leqno{(3.8.5)}
$$ 
is a resolution of~$I$ by {\it free\/} $B\b_r$-modules. 

Using the isomorphism (3.8.2), we may interpret~$\theta$ as 
the $B\b$-module homomorphism
$$
\Tor_1^{B\b}\bigl(\mu_I \bigr) \:
\Tor_1^{B\b}(I, \RR)
\longto
\Tor_1^{B\b}(I, \RR)
$$ 
induced by $\mu_I \: I \to I$, the multiplication with $T$. 
Now $\Tor_1^{B\b}(I, \RR)$ is a finite dimensional graded
vector space over~$\RR$. Hence, the injective endomorphism~$\theta$ 
has to be surjective. On the other hand, $\mu_I$ and thus 
$\theta = \Tor_1^{B\b}\bigl(\mu_I \bigr)$ has degree $2$, so it is 
not surjective unless $\Tor_1^{B\b}(I, \RR) = 0$. This yields the 
desired vanishing result~(3.8.1).
\qed

\smallskip
We still have to state and prove the two lemmata referred to above. 
The first one is a general result of commutative algebra.

\medskip
\noindent 
{\bf 3.12 Lemma:} {\it Let $R$ be a polynomial algebra over a field 
and consider an exact sequence 
$$ 
  0 \longto R^s  \longto M
  \longto L \longto 0
$$
of $R$-modules. If~$M$ is torsion free and finitely generated, then 
either $L=0$ or the codimension of its support $\supp(L)$ in the 
spectrum ${\Spec}\,R$ is at most~$1$.}

\smallskip
\noindent
{\it Proof:} We may assume that $Y := \supp\,L$ is a proper subset of 
$X := \Spec\,R$. Hence~$L$ is a torsion module, and thus~$M$ is of 
rank~$s$. Let~$Q$ be the field of fractions of~$R$. Since~$M$ is 
torsion-free, there is a natural monomorphism 
$$
  M = M \otimes_{R} R \into M \otimes_{R} Q =: M_{Q} \cong Q^s \;.
$$
We may interpret the given monomorphism $R^s \into M$ as an inclusion. 
Hence, an $R$-basis of $R^s$ may be considered as a $Q$-basis of 
$M_{Q}$, thus providing an identification $M_{Q} = Q^s$. 

We assume $\codim_X (Y) \ge 2$ and show $l=0$ for every element
$l \in L$. So fix an inverse image $m = (q_{1}, \dots, 
q_{s}) \in M \subset Q^s$ of that element $l \in L$. 
A prime ideal~$\pp$ of~$R$ lies in 
$X \setminus Y$ if and only if the localized module~$L_{\pp}$ 
vanishes, or equivalently -- since localization is exact --, if and 
only if the localized inclusion $(R_{\pp})^s \into M_{\pp}$ is an 
isomorphism. Hence, $\pp \not\in Y$ implies $q_{1}, \dots, q_{s} 
\in R_{\pp}$. Since a polynomial ring over a field is normal, the 
stipulation $\codim_{X}(Y) \ge 2$ yields $q_{1}, \dots, q_{s} \in R$
and hence $m \in R^s$, thus proving $l = \{0\}$.\qed

\medskip
\noindent
{\bf 3.13 Lemma:} {\it The support of the $A\b$-module 
$\tilde H^q(\Delta, \partial \Delta; \E\b)$ has codimension \par
\noindent
$c \ge q+2$ 
in ${\rm Spec}(A\b)$.}

\smallskip
\noindent
{\it Proof:} We show that the support is contained in the union of the 
``linear subspaces" $\Spec(A\b_\sigma)$ of $\Spec(A\b)$ with 
$\dim \sigma \le n-q-2$. To that end, we consider a prime ideal 
$\pp \in \Spec(A\b)$. Since 
localization of $A\b$-modules at~$\pp$ is exact, the localized 
cohomology module $\tilde H^q_{\pp}$ is the $q$-th cohomology of the 
complex 
$$
\tilde C\b_{\pp} \widecong
\tilde C\b (\Delta, \partial \Delta;
\E\b_{\pp})\ ,
$$
where the ``localized'' sheaf $\E\b_{\pp}$ is defined by setting 
$$
  \E\b_{\pp} (\tau) := \E\b (\tau)_{\pp}\ .
$$
Let $k$ be the minimal dimension of a cone $\tau \in \Delta$ such 
that~$\pp$ belongs to $\Spec (A_\tau)$. Then $\E\b_{\pp}(\sigma) = 0$ 
for a cone with $\dim \sigma < k$, hence in particular
$$
\E\b_{\pp} \cong \bigoplus_{\dim \sigma \ge
k} {}_\sigma \J \otimes K_\sigma
$$ 
with the characteristic sheaves ${}_\sigma \J$ and
suitable vector spaces
$K_\sigma$ and thus, according to (3.6.1) and
Remark 3.7 
$$
\tilde H^q(\Delta, \partial \Delta;
\E\b_{\pp})
\cong
\bigoplus_{\dim \sigma \ge
k} \tilde H^q (\Delta, \partial
\Delta; {}_\sigma \J) \otimes K_\sigma =0\
{\rm for}\  q > n-k-2\ .
$$
Assume now $\tilde H^q(\Delta, \partial \Delta;
\E\b)_{\pp} \cong \tilde H^q(\Delta, \partial \Delta;
\E\b_{\pp}) \not= \{0\}$ for some $\pp$ not contained in the union
of the linear subspaces $\Spec (A\b_\sigma)$ with $\dim \sigma
\le n-q-2$. Thus we have $k >n-q-2$ resp. $q >n-k-2$, a contradiction.
So $\supp(\tilde H^q(\Delta, \partial
\Delta; \E\b))$ is contained in the union
of the ``linear subspaces"
$\Spec(A\b_\sigma)$ with $\dim \sigma
\le n-q-2$, as was to be proved. 
\qed

\medskip
Eventually, we come back to the relation with the minimal complexes 
in the sense of Bernstein and Lunts: In [BeLu], a complex 
$$
  K\b : 0 \longto K^{-n} \buildrel{\delta^{-n}}\over\longto K^{-n+1} 
  \buildrel{\delta^{-n+1}}\over\longto \dots 
  \buildrel{\delta^{-1}}\over\longto K^0 \longto 0 
$$
of graded $A\b$-modules is called {\it minimal\/} 
if it satisfies the following conditions:
\item
{(i)} $K^0 \cong \RR\b[n]$, i.e., the $A\b$-module $A\b/\mm \cong 
      \RR\b$ placed in degree $-n$;
\item
{(ii)} there is a decomposition $K^{-d} = 
       \bigoplus_{\sigma \in \Delta^d} K_{\sigma}$ for $0 \le d \le n$;
\item
{(iii)} each $K_{\sigma}$ is a free graded $A\b_{\sigma}$-module; 
\item
{(iv)} for each cone $\sigma \in \Delta$, the differential $\delta$ 
       maps $K_{\sigma}$ to $\bigoplus_{\tau \prec_{1} \sigma}K_{\tau}$, 
       so for $\dim\sigma = d$, one obtains a subcomplex  
$$
  0 \longto K_{\sigma} \buildrel{\delta^{-d}_{\sigma}}\over\longto 
  \bigoplus_{\tau \prec_{1}\sigma}K_{\tau} 
  \buildrel{\delta^{-d+1}_{\sigma}}\over\longto \dots \longto K_{o} 
  \longto 0 \;;
$$ 
\item
{(v)} with $I_{\sigma} := \ker\,\delta^{-d+1}_{\sigma}$, 
      the differential  $\delta^{-d}_{\sigma}$ induces an isomorphism 
$$
  \quer\delta^{\,-d}_{\sigma} : 
  \quer K_{\sigma} := K_{\sigma}/\mm
  K_{\sigma} 
  \buildrel\cong \over\longto
  \quer I_{\sigma} := 
   I_{\sigma}/\mm I_{\sigma}
$$
of real vector spaces.

If the fan~$\Delta$ is purely $n$-dimensional, then the shifted 
cochain complex $K\b := C\b(\Delta, \E\b[n])[n]$ -- i.e., given by 
$K^{-i} := C^{n-i}(\Delta, \E\b[n])$ -- is minimal: 
With $K_{\sigma} := E\b_{\sigma}[n]$, conditions (i) -- (iv) 
are immediate; condition (v) follows from (LME) using the isomorphism 
$I_{\sigma} \cong \E\b (\partial \sigma)[n]= E\b_{\partial\sigma}[n]$
of $A\b_{\sigma}$-modules.

\smallskip
Theorem~3.8 provides a characterization of quasi-convex fans in terms of 
acyclicity of the relative cellular cochain complex. An analoguous statement 
holds also for the absolute cellular cochain complex. In particular, this 
proves a conjecture of Bernstein and Lunts (see [BL], p.129, 15.9):

\medskip \noindent    
{\bf 3.14 Theorem:} {\it A purely $n$-dimensional fan~$\Delta$ is 
quasi-convex if and only if the complex $C\b(\Delta, \E\b)$ is exact in degrees
$q>0$ and $H^0(\Delta,\E\b) \cong E\b_{(\Delta, \partial \Delta)}$. 
Specifically, for a complete fan~$\Delta$, a minimal complex in the 
sense of Bernstein and Lunts is exact.}  

\smallskip
\noindent
{\it Proof:} We consider the augmented absolute cellular cochain
complex
$$
0 \longto
F_{(\Delta, \partial \Delta)}
\longto 
C^0 (\Delta; \F )
\longto \dots
\longto 
C^n (\Delta; \F) 
\longto 0\ \leqno(3.14.1)
$$ 
for some sheaf $\F$ on $\Delta$. By~3.6.1, it is acyclic for the flabby 
sheaf $\E\b$ if and only if it is acyclic for each characteristic sheaf  
${}_\sigma \J$, where $\sigma \in \Delta$. For $\sigma \not\in \partial 
\Delta$, that follows from 3.5.1 and Rem.3.5, since 
$\Delta_\sigma$ is complete.
For a cone $\sigma \in 
\partial\Delta$, the absolute versions of Remark~3.7 and formula 
(3.10.0) yield isomorphisms 
$$
\tilde H^q(\Delta, {}_\sigma \J)
\cong
H^q(\Delta,{}_\sigma \J)
\cong
H^q(\Delta_\sigma, \RR)
\cong
\tilde H_{k-1-q}(L_\sigma, \RR)\ ,
$$
where $n-k= \dim \sigma$ and $L_\sigma$ is the link of some point 
$x \in S_\Delta \cap \Rond{\sigma}$.  Now statemant iii) of Proposition~3.10
gives $\tilde H_{k-1-q}(L_\sigma, \RR)= \{0\}$. \qed

\smallskip
For later use we still need the following result.

\medskip
\noindent 
{\bf 3.15 Corollary.} {\it For a minimal extension sheaf $\E\b$ on a 
quasi-convex fan $\Delta$, the $A\b$-submodule $E\b_{(\Delta, 
\partial\Delta)} \subset E\b_\Delta$ of global sections vanishing on the 
boundary fan $\partial \Delta$ is free.}

\smallskip
\noindent
{\it Proof:}
From the acyclicity of the absolute cellular cochain complex,  
we conclude as in the proof of Theorem 3.8 that 
$E\b_{(\Delta, \partial \Delta)}$ is a free $A\b$-module. \qed

\bigskip\medskip\goodbreak
\centerline{\XIIbf 4. Poincar\'e Polynomials}
\bigskip\nobreak\noindent
In the remaining part of our article, we
want to discuss the virtual intersection Betti numbers 
$b_{2q}(\Delta) :=  \dim \quer
E^{2q}_\Delta$ and 
$b_{2q}(\Delta, \partial\Delta) := \dim
\quer E^{2q}_{(\Delta, \partial \Delta)}$
of a quasi-convex  fan~$\Delta$, where
$\E\b$ is a minimal extension sheaf
on~$\Delta$.  It is convenient to use the
language of Poincar\'e polynomials. 

\bigskip \noindent
{\bf 4.1 Definition:} {\it The (
equivariant) Poincar\'e  series of a
fan~$\Delta$ is the formal power series 
$$ 
 Q_{\Delta} (t) := \sum_{q \ge
0} \dim\,
  E_{\Delta}^{2q} 
  \cdot t^{2q}\; ,
$$
its (intersection) Poincar\'e polynomial is the polynomial
$$
  P_{\Delta} (t) := \sum_{q \ge 0}^{<\infty} \dim\,
 \quer E_{\Delta}^{2q} 
  \cdot t^{2q} = \sum_{q \ge 0}^{<\infty} b_{2q}(\Delta) t^{2q}\; .
$$
For an affine fan $\langle \sigma \rangle$, we
simply write
$$
Q_\sigma := Q_{\langle \sigma
\rangle}\ ,\ 
P_\sigma := P_{\langle \sigma
\rangle} \ .
$$
Furthermore, for a subfan $\Lambda \preceq \Delta$, 
the relative Poincar\'e polynomial $P_{(\Delta, \Lambda)}$
is defined in an analoguous manner.}
\par

\bigskip 
We refer to $P_\Delta$ as the {\it global\/} Poincar\'e polynomial 
of~$\Delta$, while the polynomials $P_\sigma$ for $\sigma \in \Delta$ 
are called its {\it local\/} Poincar\'e polynomials. 
\par

\bigskip
\noindent
{\bf 4.2 Remark.} {\it For a quasi-convex fan we
have
$$
 Q_{\Delta}(t) = {1 \over \left(1-t^2\right)^n} \mal P_{\Delta}(t) 
  \ ,
$$
while for a cone $\sigma$, one has
$$
Q_{\sigma}(t) = 
  {1 \over (1-t^2)^{\dim \sigma}} \mal
  P_{\sigma}(t)\ .
$$}

\bigskip
\noindent
{\bf Proof.} For a free graded
$A\b$-module $F\b$, the K\"unneth
formula $F\b \cong A\b \otimes_\RR
\quer F\b$ holds, while the Poincar\'e series of a 
tensor product of graded 
vector spaces is the product of the
Poincar\'e series of the factors.  Since
$Q_{A\b} = 1/(1-t^2)^n$, the first
formula follows immediately.  Going over
to the base ring~$A\b_{\sigma}$
yields the second one. \qed

\bigskip
The basic idea for the computation of the virtual intersection Betti 
numbers is to use a two-step procedure. In the first step, the 
global invariant is expressed as a sum of local terms. In the second 
step, these local invariants are expressed in terms of the global 
ones associated to lower-dimensional 
fans.

\bigskip
\noindent
{\bf 4.3 Theorem (Local-to-Global
Formula):} If $\Delta$ is a
quasi-convex fan of dimension~$n$ and $\buildrel \circ \over
\Delta:=\Delta \setminus \partial \Delta$,,
we have  
$$
  P_{\Delta}(t) = \sum_{\sigma \in 
\buildrel \circ \over
\Delta} 
  (t^2-1)^{n- \dim \sigma} P_{\sigma}(t) 
  \ ,
$$
while
$$
P_{(\Delta, \partial \Delta)}(t)
=\sum_{\sigma \in 
\Delta} 
  (t^2-1)^{n- \dim \sigma} P_{\sigma}(t) 
  \ .
$$

\bigskip
\noindent
{\it Proof.} The cellular cochain complex  
$$
0 \longto E\b_\Delta \longto
C^0(\Delta, \partial \Delta ; \E\b)
\longto \dots  \longto C^n(\Delta, \partial
\Delta; \E\b) \longto 0 
$$ 
of~3.4 associated to the quasi-convex fan~$\Delta$ is acyclic by 
Theorem~3.8. We set 
$$
Q_i(t) := 
\sum_{q \ge 0} \dim\,
  C^i(\Delta,\partial \Delta ; \E^{2q}) 
  \cdot t^{2q}=\sum_{\sigma \in 
{\buildrel \circ \over \Delta} \cap
\Delta^{n-i}} Q_\sigma (t) \;.
$$
Then we obtain the equality 
$$
Q_\Delta (t)=\sum_{i=0}^n
(-1)^iQ_i(t)=\sum_{\sigma \in
\buildrel \circ \over \Delta}
              (-1)^{n - \dim \sigma}
              Q_\sigma (t)\ .
$$
The first assertion follows from Remark 4.2. The second formula 
is obtained in the same way using the acyclicity of the complex
$$
0 \longto E\b_{(\Delta, \partial \Delta)}
\longto C^0(\Delta,
\E\b) \longto \dots  \longto C^n(\Delta,
\E\b) \longto 0 
$$
see Theorem~3.14 and Corollary~3.15. \qed 

\bigbreak
In order to reduce the computation of $\quer E\b_\sigma \cong
\quer
E\b_{\partial \sigma}$ to a problem in lower dimensions, we choose a
line $L \subset V$ meetimg the relative interior $\buildrel \circ
\over \sigma$ and consider the flattened boundary
fan $\Lambda_\sigma:= \pi (\partial
\sigma)$, where $\pi: V_\sigma \to V_\sigma /L$ is the
quotient map, cf. section 0.D. Then the direct image sheaf
$$
  \G\b := \pi_{*}(\E\b|_{\partial\sigma}) :\; \tau \mapsto 
  \E\b\big((\pi|_{\partial\sigma})^{-1}(\tau)\big) \leqno{\rm(4.3.1)}
$$
for $\tau \Lambda_{\sigma}$, is a minimal
extension sheaf on $\Lambda_\sigma$. 
Writing $A\b_\sigma = B\b[T]$ with 
$B\b:=\pi^*(S\b((V/L)^*)) \subset A\b_\sigma$ and
$T \in A^2_\sigma$ as in section~0.D, we 
obtain the identification 
$$
  \quer E\b_{\sigma} \cong \quer E\b_{\partial\sigma} \cong 
  \quer G\b_{\Lambda_{\sigma}}/(f \!\cdot\! \quer 
  G\b_{\Lambda_{\sigma}})\leqno{\rm(4.3.2)}
$$
with the piecewise linear function $f:=T \circ 
(\pi|_{\partial\sigma})^{-1} \in \A^2(\Lambda_{\sigma})$.
Here $E\b_\sigma$ and $E\b_{\partial \sigma}$ are considered as
$A\b_\sigma$-modules, while $G\b_{\Lambda_\sigma}$ is a $B\b$-module
only; and it is with respect to that module structures one has to take
residue class vector spaces.

As a first result we get an estimate for the degree of the Poincar\'e 
polynomials:

\bigskip
\noindent
{\bf 4.4 Corollary:} \hangafter=1\hang {\it 
i) For a quasi-convex fan
$\Delta$ the relative
Poincar\'e polynomial $P_{(\Delta,
\partial \Delta)}$ is  monic of
degree~$2n$, whereas for a non-complete
quasi-convex fan $\Delta$ the absolute
Poincar\'e polynomial $P_\Delta$ 
is of
degree at most $2n-2$.
\item{ii)} For a non-zero cone~$\sigma$,
the ``local'' Poincar\'e  polynomial
$P_\sigma$ is of degree at most $2 \dim
\sigma - 2$.}

\bigskip
\noindent
{\it Proof.} We proceed by induction on the dimension. Assuming that~(ii) 
holds for every cone~$\sigma$ with $\dim\sigma \le n$, then Theorem~4.3 
yields~(i) in dimension~$\dim\Delta = n$. Now, if $\sigma$ is a ray, 
assertion~(ii) is evident. For the induction step, we now assume 
$\dim \sigma = n > 1$. Going over to the complete 
fan~$\Lambda_{\sigma}$ of dimension~$n-1$, we use the isomorphism 
(4.3.2). Since $\quer G^q_{\Lambda_{\sigma}} = 0$ holds for $q > 
2n-2$ according to the induction
hypothesis, assertion~(ii) follows.\qed

\bigskip
For the second step, we have to relate the local Poincar\'e polynomial 
$P_\sigma$ to the global Poincar\'e polynomial 
$P_{\Lambda_\sigma}(t)$ of the complete (and thus quasi-convex) fan 
$\Lambda_\sigma$ of dimension $\dim \sigma - 1$.
Here the vanishing condition
$V(\sigma)$, cf.1.7, plays a decisive
role:
\par

\bigskip
\noindent
{\bf 4.5 Theorem (Local Recursion Formula):} {\it Let $\sigma$ be a cone. 
\item{i)} If $\sigma$ is simplicial, then we have $P_\sigma \equiv 1$. 
\item{ii)} If the condition $V(\sigma)$ is satisfied and $\sigma$ is 
not the zero cone, then we have 
$$ 
  P_{\sigma}(t) = 
  \tau_{<\dim\sigma}\big((1-t^2)
  P_{\Lambda_\sigma}(t)\big)\;.
$$
}

In the statement above, the truncation operator $\tau_{<k}$
is defined by $\tau_{<k}(\sum_q a_q t^q):=\sum_{q < k} a_q t^q$. -- 
Let us note that for $\dim\sigma = 1$ and $2$, the statements (i) and 
(ii) agree. 

\bigskip
\noindent
{\it Proof:} Statement i) follows from
the fact that $E\b_\sigma \cong
A\b_\sigma$ for a simplicial cone $\sigma$. In order to prove
statement ii) we use the isomorphism (4.3.2), thus have to investigate
the graded vector space $\quer G\b_{\Lambda_\sigma}/f
\quer G\b_{\Lambda_\sigma}$ respectively the kernel and cokernel
of the map
$$
\quer \mu: \quer G\b_{\Lambda_\sigma}[-2]
\longto \quer G\b_{\Lambda_\sigma}, \quer h \mapsto
\quer {fh}
$$
induced by the multiplication $\mu \: G\b_{\Lambda_\sigma}[-2]
\to G\b_{\Lambda_\sigma}$ with the piecewise linear function $f \in 
\A^2(\Lambda_\sigma)$. We apply the following ``Hard Lefschetz'' 
type theorem with $\Delta=\Lambda_\sigma$ and $\E\b = \G\b$:

\bigskip 
\noindent
{\bf 4.6 Theorem:} {\it Let $\Delta$ be a
complete fan and $f \in \A^2(\Delta)$
a strictly convex function, i.e.,
$f=(f_\sigma)_{\sigma \in \Delta}$
is convex and $f_\sigma \not=
f_{\sigma'}$ for $\sigma \not=
\sigma'$, furthermore let $\gamma^+(f) \subset V \times \RR$ 
be the convex hull of
the graph $\Gamma_f \subset V \times \RR$ of $f$. Then, if the condition  
$V\big(\gamma^+(f))$ is satisfied, the map 
$$ 
\quer \mu^{2q}:\quer 
E^{2q}_{\Delta}
\longto \quer
E^{2q+2}_{\Delta} 
$$
induced by the multiplication 
$\mu: E\b_{\Delta}[-2] \to
E\b_{\Delta} \ , \ h \mapsto  fh
$, is injective for $2q \le n-1$ and
surjective for $2q \ge n-1$.}

\bigskip
Theorem 4.6 will be derived from the vanishing condition $V(\gamma^+(f))$
at the end of 
section~5 by means of  
Poincar\'e duality. As a simple consequence 
of Corollary~5.6, a ``numerical''
version of Poincar\'e duality can be formulated as 
follows.

\bigskip
\noindent
{\bf 4.7 Theorem:} {\it For a quasi-convex fan~$\Delta$, 
the global Poincar\'e polynomials $P_\Delta$ and 
$P_{(\Delta,\partial\Delta)}$ satisfy the identity
$$ 
P_{(\Delta,\partial\Delta)}(t) = t^{2n} P_\Delta(t^{-1})\; . 
$$}
\bigskip
We conclude this section with an
application of the decomposition theorem
2.3, which has been communicated to us by
Tom Braden (cf. also [BrMPh]):

\medskip
\noindent
{\bf 4.8 Theorem (Kalai's
conjecture)} {\it For a face $\tau \preceq
\sigma$ of the cone $\sigma$ with
transversal fan $\Delta_\tau$ in
$\Delta:=\langle \sigma \rangle$ we
have:   
$$ 
P_\sigma (t) \ge P_\tau (t)
\cdot P_{\Delta_\tau} (t) \ ,  
$$ 
where $P \ge Q$
means, that the corresponding inequality
holds for the coefficients of monomials of $P$ and
$Q$ with the same degree.}
\bigskip

\noindent
{\it Proof.} Consider the
minimal extension sheaf $\E\b$ on the
affine fan $\Delta:=\langle \sigma
\rangle$ and denote $\F\b$ the trivial
extension of $\E\b|_{{\rm st}(\tau)}$,
i.e. if $\Lambda \preceq \Delta$ and
$\Lambda_0$ is generated by the cones in
$\Lambda \cap {\rm st}(\tau)$, then $\F\b
(\Lambda) = \E\b (\Lambda_0)$. Obviously
$\F\b$ is a pure sheaf and its
decomposition has the form $$\F\b
\cong ({}_\tau \E\b \otimes
\overline  E\b_\tau)
\oplus \bigoplus_{\gamma \succ \tau} 
{}_\gamma \E\b \otimes 
K\b_\gamma\ .
$$
Now our inequality follows immediately
by taking the residue class module of the
global sections of the above sheaf, since
${}_\tau\E\b \cong 
B\b_\tau \otimes_\RR
({}_{\Delta_\tau} \E\b)$ (identifying
$\Delta_\tau$ with st$(\tau)$) with
$B\b:=S((V/V_\tau)^*) \subset A\b$ and
$B\b_\tau \subset A\b$, such that
$B\b_\tau \cong A\b_\tau$, hence in
particular $A\b \cong B\b_\tau \otimes_\RR
B\b$.  \qed
\bigskip

\bigskip\medskip\goodbreak
\centerline{\XIIbf 5. Poincar\'e Duality}
\bigskip\nobreak\noindent
In this section, we first define a -- non-canonical -- ``intersection 
product'' $\E\b \times \E\b \to \E\b$ on a minimal extension sheaf 
$\E\b$ for an arbitrary fan~$\Delta$. On the level of 
global sections, it provides a ``product'' $E\b_{\Delta} \times 
E\b_{(\Delta, \partial \Delta)} \to E\b_{(\Delta, \partial \Delta)}$. 
If the fan is even 
quasi-convex, then in addition, there exists an evaluation mapping 
$\epsilon: E\b_{(\Delta, \partial\Delta)} \to A\b[-2n]$. The main 
result of this section is the ``Poincar\'e Duality Theorem''~5.3 
according to which the composition of the intersection product 
and the evaluation map is a dual pairing. 

In the case of a simplicial fan, where the sheaf~$\A\b$ of piecewise 
polynomial functions is a minimal extension sheaf, such a product is 
simply given by the multiplication of these functions. Hence, one 
possible approach to the general case is as follows: Choose a 
simplicial refinement $\hat \Delta$ of $\Delta$ and, according to the 
Decomposition Theorem~2.4,
interpret~$\E\b$ as a direct factor of
the sheaf $\hat \A$ of $\hat
\Delta$-piecewise polynomial functions 
on~$\Delta$. Then take the restriction of
the multiplication of functions  on $\hat
\A\b$ to $\E\b \subset \hat \A\b$ and
project to~$\E\b$.  \par

But in order to keep track of the relation
between the intersection  product over
the boundary of a cone and the cone
itself, it is useful  to  apply the above
idea repeatedly in a recursive
extension procedure. The proof of
Poincar\'e  duality will follow the same
pattern.

\bigskip 
\noindent 
{\bf 5.1 An Intersection Product:} The $2$-dimensional skeleton 
$\Delta^{\le 2}$ is a simplicial subfan. Hence, up to a scalar multiple, 
there is a canonical isomorphism $\A\b \cong \E\b$ on $\Delta^{\le 2}$ 
(see 1.8). We thus define the intersection product on $\Delta^{\le 2}$ 
to correspond via that isomorphism to the product of functions.
\par

We now assume that the intersection product is defined on $\Delta^{\le m}$
and consider a cone $\sigma \in \Delta^{m+1}$. So we are given a symmetric 
bilinear morphism $E\b_{\partial\sigma} \times E\b_{\partial\sigma} \to 
E\b_{\partial\sigma}$ of $A\b_{\sigma}$-modules. As in the previous
section, let $L \subset V_\sigma$ be a line intersecting $\buildrel \circ
\over \sigma$ and let $B\b \subset A\b_\sigma$ be the image of 
$S\b((V_\sigma/L)^*)$ in $A\b_\sigma = S(V_\sigma^*)$. We recall that 
$E\b_{\partial \sigma} \cong G\b_{\Lambda_\sigma}$ (see the remarks
preceding Corollary~4.4) is a free $B\b$-module, by Theorem~3.8. Let us 
define the $A\b_\sigma$-module 
$$
  F\b_\sigma \;:=\; A\b_\sigma \otimes_{B\b} E\b_{\partial \sigma}\;. 
  \leqno(5.1.1)
$$
Since $E\b_\sigma$ is a free $A\b_\sigma$
-module,  
the (surjective) 
restriction $E\b_\sigma \to E\b_{\partial \sigma}$ can be factorized
$$
  E\b_{\sigma} \buildrel \alpha \over \longto 
  A\b_\sigma \otimes_{B\b} E\b_{\partial\sigma} = 
  F\b_\sigma  \buildrel \beta \over \longto 
  A\b_\sigma \otimes_{A\b_\sigma} E\b_{\partial\sigma} = 
  E\b_{\partial\sigma} \;.
$$
The map $\alpha \: E\b_{\sigma}
\to F\b_\sigma$ is a ``direct''
embedding,  i.e., there is a decomposition
$F\b_\sigma \cong \alpha (E\b_{\sigma})
\oplus K\b$, since the reduction
of~$\alpha$ mod~$\mm \subset A\b_\sigma$ 
is injective and $F\b_\sigma$ a free
$A\b_\sigma$-module. We may even  assume
that $K\b$ is contained in the kernel of
the natural map  $A\b_\sigma
\otimes_{B\b} E\b_{\partial \sigma} \to
E\b_{\partial \sigma}$:  Take a
homogeneous basis $f_1,\dots ,f_r$ of
$K\b$. The images $\beta(f_{i})$ of these
elements in $E\b_{\partial \sigma}$ are
restrictions of elements $g_i \in
E\b_\sigma$; hence, we may  replace $K\b$
by the submodule generated by the
elements $f_i-\alpha (g_i)$ for $1 \le i
\le r$. \par

On the other hand, by scalar extension, there is an induced product 
$$
F\b_\sigma
\times
F\b_\sigma
\longto
F\b_\sigma\ . 
$$
It provides the desired extension of the intersection product 
from $\partial \sigma$ to $\sigma$ 
via the 
composition 
$$
E\b_\sigma \times E\b_\sigma
\;\buildrel \alpha \times \alpha \over
\longto\; F\b_\sigma \times
F\b_\sigma \longto
F\b_\sigma = \alpha
(E\b_\sigma) \oplus K\b  
\longto \alpha (E\b_\sigma)
\cong E\b_\sigma\ ,
$$
where the last arrow is the projection
onto $\alpha (E\b_\sigma)$ with kernel
$K\b$. This ends the extension procedure. To sum up, after a finite 
number of steps, we arrive at a symmetric bilinear morphism 
$$
\E\b \times \E\b \longto \E\b \leqno(5.1.2)
$$
of sheaves of $\A\b$-modules, called an {\it intersection product
on the minimal extension sheaf $\E\b$\/}. In particular, we thus have 
defined a product 
$$
E\b_\Delta \times E\b_\Delta 
\longto E\b_\Delta\; ,
$$
and obviously we have $E\b_\Delta \cdot E\b_{(\Delta,
\partial \Delta)} \subset 
E\b_{(\Delta,
\partial \Delta)}$.

\medskip
In order to obtain a dual pairing in the
case of a quasi-convex fan~$\Delta$, we
compose the induced product $E\b_\Delta \times
E\b_{(\Delta, \partial \Delta)} \to
E\b_{(\Delta,
\partial \Delta)}$
with an
``evaluation''  homomorphism  
$$
\epsilon : E\b_{(\Delta, \partial \Delta)}
\longto A\b[-2n] 
$$
that can be defined as follows: As a consequence of Corollary~4.4, we know 
$$
\quer E^q_{(\Delta, \partial \Delta)}
= \cases{ \RR &
, $q=2n$ \cr
0 &, $q > 2n$ \cr}\; .
$$
Moreover, according to
Corollary~3.15, $E\b_{(\Delta, \partial
\Delta)}$ is a free  $A\b$-module. Hence,
there is an element  $\epsilon \in
\Hom_{A\b}(E\b_{(\Delta,\partial
\Delta)}, A\b [-2n]) \setminus \{0\}$  of
degree zero; in fact it is unique up to 
multiplication by a real scalar. If
$\Delta$ is a simplicial fan, this 
homomorphism~$\epsilon$ can be  described
quite explicitly: Following [Bri,
p.13], one  fixes some euclidean norm on
$V$ and hence also on $V^*$. Write each cone $\sigma \in
\Delta^n$ as $\sigma=H_1 \cap \dots \cap
H_n$ with half spaces $H_i=H_{\alpha_i}$
and linear forms $\alpha_i \in V^*$ of length $1$, finally
set  $f_\sigma := \alpha_1 \cdot \dots
\cdot \alpha_n$. Then the map~$\epsilon$
is of the following form:  
$$ 
  E\b_{(\Delta, \partial \Delta)} \cong
A\b_{(\Delta, \partial \Delta)} \subset
\bigoplus_{\sigma \in \Delta^n}
A\b_\sigma \longto A\b[-2n]\;,\;\;
h=(h_\sigma)_{\sigma \in \Delta^n} \mapsto
\sum_{\sigma \in \Delta^n} {h_\sigma \over
f_\sigma}\ .
$$

\smallskip
\noindent 
For example, if $\Delta = \langle\sigma\rangle$ is a full-dimensional 
affine simplicial fan, then $A\b_{(\sigma, \partial\sigma)}$ is of the 
form $f_\sigma A\b_\sigma$ for the function $f_{\sigma} 
\in A^{2n}$ as above, and hence the
above map has values in $A\b \subset
Q(A\b)$. In the general case for each
$n-1$-cone in $\buildrel \circ \over
\Delta$ there are two summands which
have a pole of order $1$ along it, but
these poles cancel one another, while
each summand already is regular along
$\partial \Delta$ because of $h \in
A\b_{(\Delta, \partial \Delta)}$. 

Since the intersection product $\E\b
\times \E\b \to \E\b$ is a  homomorphism
of sheaves, we may sum up the general
situation as  follows: For a quasi-convex
fan~$\Delta$, there exists {\it
homogenous  pairings\/} (i.e., of degree
zero with respect to the total grading
on  the product) 
$$ 
E\b_\Delta \times 
E\b_{(\Delta, \partial \Delta)}
\longto
E\b_{(\Delta, \partial \Delta)}
\longto
A\b[-2n] \leqno(5.1.3)
$$
and
$$
\quer E\b_\Delta
\times 
\quer E\b_{(\Delta, \partial \Delta)}
\longto
\quer E\b_{(\Delta, \partial \Delta)}
\longto
\RR\b[-2n]\ . \leqno(5.1.4)
$$

\bigskip
Our aim is to prove that these are in fact both dual pairings. 
Fortunately, it suffices to verify that property for one of these two 
pairings: According to the very definition of a minimal extension 
sheaf and by Theorem~3.8 and
Corollary~3.15, the $A\b$-modules 
$E\b_\Delta$ and $E\b_{(\Delta, \partial
\Delta)}$ are free. We may  thus apply
the following result. 

\bigskip
\noindent
{\bf 5.2 Lemma.} {\it Let $E\b$, $F\b$ be two 
finitely generated free graded 
$A\b$-modules. Then a homogeneous
pairing 
$$
E\b \times F\b \to A\b[r]
$$
is dual if and only if the
induced pairing
$$ 
\quer E\b \times \quer F\b
\longto \quer
A\b[r] 
$$
is.}

\medskip
\noindent
{\it Proof.} After shifting the grading of $F\b$, we may assume $r=0$. 
We choose homogeneous bases of $E\b$ and $F\b$, so the pairing can be 
represented by a matrix~$M$ over~$A\b$. Then~$M$ is a square matrix 
and is invertible if and only if that holds for its residue class 
mod~$\mm_{A}$: The implication ``$\Longrightarrow$" is obvious, while 
for ``$\Longleftarrow$", it suffices to
prove $\det M \in A^0 = \RR$.

We arrange the bases in increasing order for $E\b$ and decreasing 
order for $F\b$ with respect to the degrees. Since 
the induced pairing is a dual one,
the submodules of $E\b$
and $F\b$ generated by basis elements of opposite degrees have the 
same rank. The matrix of the pairing is composed of square matrices 
with entries in $A^0$ along the diagonal, and below these 
all entries are $0$. Thus the determinant $\det M$ 
equals the product of the determinants
of the diagonal square blocks, so it is a
constant. 
\qed

\medskip
We come now to the central result of this section: 

\medskip
\noindent
{\bf 5.3 Theorem (Poincar\'e Duality (PD)):}  {\it For a quasi-convex fan 
$\Delta$ of dimension $n$, the composition 
$$
E\b_\Delta \times 
E\b_{(\Delta, \partial \Delta)}
\longto
E\b_{(\Delta, \partial \Delta)}
\longto A\b[-2n]
$$
is a dual pairing of finitely generated
free $A\b$-modules.}

\smallskip
\noindent
{\it Proof:} For an affine fan $\Delta:=\langle \sigma \rangle$ with a 
cone of dimension $\dim \sigma =n\le 2$ Poincar\'e duality obviously 
holds. Now Theorem~5.3 follows by the following two lemmata, since using them
one can prove the general case with a two step induction procedure.
\medskip
\noindent
{\bf 5.3a Lemma.} {\it If Poincar\'e duality holds for complete fans in 
dimensions
$<n$, then also for an affine fan $\Delta =\langle \sigma \rangle$ with a cone
$\sigma$ of dimension $n$.}
\medskip
\noindent
{\it Proof.} According to (4.3.1), we identify 
$E\b_{\partial \sigma}$ with the 
$B\b$-module $G\b_{\Lambda_{\sigma}}$ of global sections 
of a minimal extension sheaf~$\G\b$ on the fan
$\Lambda_\sigma$ in $V/L$. Since the fan~$\Lambda_\sigma$ is of 
dimension $<n$,
we obtain a dual pairing
$$
E\b_{\partial \sigma}
\times E\b_{\partial \sigma}
\longto E\b_{\partial \sigma}
\longto B\b[2-2n]\;.
$$ 
By extension of scalars, that induces dual pairings 
$$
F\b_\sigma \times F\b_\sigma
\longto F\b_\sigma
\buildrel
\eta
\over
\longto A\b[2-2n]\ 
$$
resp.
$$
\quer F\b_\sigma
\times \quer F\b_\sigma
\longto \quer F\b_\sigma
\longto
\RR\b[2-2n]\ 
$$ 
and, after a shift, 
$$\
\quer F\b_\sigma \times 
\quer F\b_\sigma [-2]
\longto 
\quer F\b_\sigma [-2] \longto
\RR\b[-2n]\ . 
$$
To achieve the proof, we are going to construct a factorization of 
the induced pairing
$\quer E\b_\sigma \times
\quer E\b_{(\sigma,\partial \sigma)}
\to \RR\b[-2n]$
on the level of residue class vector spaces in the following form:  
$$
\quer E\b_\sigma \times
\quer E\b_{(\sigma,\partial \sigma)}
\buildrel {\alpha \times \theta} \over \longto
\quer F\b_\sigma \times 
\quer F\b_\sigma [-2]
\longto 
\quer F\b_\sigma [-2] \longto
\RR\b[-2n]\,. \leqno(5.3.1)
$$
We show the existence of a homomorphism~$\mu: \quer F\b_\sigma[-2]
\to \quer F\b_\sigma$ such
that~$\alpha$ and~$\theta$ induce 
isomorphisms  $$
\quer E\b_{(\sigma,\partial \sigma)}
\cong \ker\,\mu  \qquad\hbox{and} \qquad
\quer E\b_\sigma
\cong \coker \,\mu \; .
$$
Finally, forgetting about the shifts, the map $\mu$ is self-adjoint 
with respect to the above dual pairing on $\quer F\b_\sigma$. 
Hence, the restriction of the pairing to 
$\ker\, \mu \times \coker\, \mu$ is dual, too. 
\par

We interpret $F\b_\sigma$ as 
the module of sections of a sheaf of
$\A\b$-modules on the affine fan $\langle
\sigma \rangle$. To that end, we consider the subdivision 
$$
\Sigma := \partial \sigma \cup \{
\hat \tau := \tau + \varrho; \tau \in
\partial \sigma \}\ .
$$
of $\langle \sigma \rangle$, where~$\rho$ is the ray $L \cap \sigma$. 
Let $B\b_\tau \subset A\b_{\hat \tau}$ denote the subalgebra of 
functions constant on parallels to the line $L$.
Then, according to remark 1.5, the sheaf
$\F\b$ on $\Sigma$ with $$
\tau \mapsto F\b_\tau:=E\b_\tau \,,\;
\hat \tau \mapsto 
F\b_{\hat \tau}:=
A\b_{\hat \tau}
\otimes_{B\b_\tau} E\b_\tau \qquad\hbox{for $\tau \in \partial \sigma$}
$$
and the obvious restriction homomorphisms is a minimal
extension sheaf $\F\b$ on $\Sigma$ and satisfies 
$\F\b(\Sigma) \cong
A\b \otimes_{B\b} E\b_{\partial \sigma} = F\b_\sigma$. 
Furthermore, 
the sheaf~$\F\b$ inherits an intersection product 
from~$\E\b|_{\partial\sigma}
\cong \F\b|_{\partial\sigma}$ as in 5.1.
\par

For simplicity, we interpret the mapping~$\alpha$ in~5.1 as an 
inclusion $E\b_{\sigma} \subset F\b_{\sigma}$ and 
identify $\F\b$ 
with its direct image sheaf on the affine
fan $\Delta:=\langle 
\sigma \rangle$ (with respect to the refinement mapping
$\Sigma \to \langle \sigma \rangle$).
Then the decomposition
$F\b_\sigma = E\b_\sigma 
\oplus K\b$ corresponds to a
decomposition $\F\b \cong \E\b \oplus
\K\b$ with $\E\b \cong {}_o \E\b$ and
the ``skyscraper'' sheaf $\K\b:=
{}_\sigma \E\b \otimes K\b$
supported by $\{ \sigma \}$. In
particular, there is an inclusion   
$$
E\b_{(\sigma, \partial \sigma)}
\subset F\b_{(\sigma, \partial \sigma)}
= E\b_{(\sigma, \partial \sigma)} \oplus
K\b \;, 
$$
and $F\b_{(\sigma, \partial\sigma)}$ is a free $A\b$-module. 

We thus obtain a natural commutative diagram
$$
\diagram{0 & 
\longto & E\b_{(\sigma, \partial
\sigma)} & \longto &
E\b_{\sigma} &
\longto & 
E\b_{\partial \sigma} &
\longto & 0\cr 
 & & \cap & & \cap & & \Vert & \cr
0 & \longto & 
F\b_{(\sigma, \partial \sigma)} & 
\buildrel \lambda
\over
\longto & 
F\b_\sigma & 
\longto & 
F\b_{\partial \sigma} & 
\longto & 0 \cr}
$$
consisting of free resolutions of the $A\b$-module 
$E\b_{\partial\sigma} \cong F\b_{\partial\sigma}$. 

Using the very definition of $\Tor^{A\b}(*, \RR\b)$ and the fact that
$\quer E\b_{(\sigma, \partial \sigma)} \to \quer E\b_\sigma$ is the 
zero map since $\quer E\b_\sigma \to \quer E\b_{\partial \sigma}$ 
is an isomorphism, we obtain identifications 
$$
\quer E\b_{(\sigma,\partial \sigma)} \widecong
\Tor_1(E\b_{\partial \sigma}, \RR\b)
\widecong \ker (\quer \lambda) \quad\hbox{and}\quad
\quer E\b_\sigma
\widecong \coker (\quer \lambda) \widecong
\quer E\b_{\partial \sigma}\;. \leqno(5.3.2)
$$

On the other hand, we may rewrite $F\b_{(\sigma, \partial \sigma)} 
= g F\b_\sigma \cong F\b_\sigma[-2]$, where $g \in \A^2(\Sigma)$ is 
some piecewise linear function on $\Sigma$ with $\partial \sigma$ as 
zero set: Write $A\b = B\b[T]$
with $B\b:=S\b((V/L)^*) \subset A\b$, such
that the kernel of $T \in A^2$  meets
$\sigma$ only in~$0$. Then, for $\tau \in
\partial\sigma$, we  set $g_{\hat \tau} =
T-f_\tau$, where $f_\tau \in A^2_{\hat
\tau}=A^2$  coincides with~$T$ on~$\tau$
and is constant on parallels to~$L$,
i.e., $f_\tau \in B\b$. 
\par

Now note that 
$$
E\b_{(\sigma,
\partial \sigma)} \subset 
E\b_{(\sigma,
\partial \sigma)}
\oplus K\b =
F\b_{(\sigma,\partial \sigma)} 
=
gF\b_\sigma \cong
F\b_\sigma [-2]
\buildrel \eta [-2]
\over \longto
A\b [-2n]
$$ 
defines the homomorphism $\theta \: \quer E\b_{(\sigma, \partial\sigma)} 
\to \quer F\b_{\sigma}[-2]$ mentioned above and an evaluation map for 
$E\b_{(\sigma, \partial \sigma)}$, and 
that $K\b \subset F\b_{(\sigma, \partial\sigma)}$ is contained in the 
kernel
of the map $F\b_{(\sigma,\partial
\sigma)}  \to
A\b [-2n]$,
since $\quer K^q=0$ for $q \ge 2n$ because of the isomorphism 
$\quer E^{2n}_{(\sigma, \partial \sigma)}
 \cong \RR \cong 
\quer F^{2n}_{(\sigma,\partial \sigma)}
\cong
\quer F^{2n-2}_\sigma  
$
and the vanishing $\quer F^q_{(\sigma,\partial \sigma)}=0$ for $q > 2n$.
Next we remark that, although the first part of the diagram 
$$
\diagram{
E\b_\sigma \times
E\b_{(\sigma,\partial \sigma)} &
\longto &
E\b_{(\sigma,\partial \sigma)} & 
\buildrel \epsilon \over \longto &
A\b [-2n] \cr
\cap && \cap && \Vert \cr
F\b_\sigma \times
F\b_{(\sigma,\partial \sigma)} &
\longto &
F\b_{(\sigma,\partial \sigma)} & 
\buildrel \eta [-2] \over \longto &
A\b [-2n] \cr}
$$
need not be commutative ($E\b_\sigma$ is not necessarily 
closed under the intersection product in
$F\b_\sigma$), commutativity holds after 
evaluation (where the two evaluation maps
are scaled in such a way  that the right
square is commutative). This is true 
since the difference of the products in
the first and second row is an element in
$K\b$, according to the construction. \par

As the intersection product in 
$F\b_\sigma$ is $\A\b(\Sigma)$-linear,
we may replace $F\b_{(\sigma,
\partial \sigma)}$ with $F\b_\sigma [-2]$
and arrive at the following pairing of $A\b$-modules
$$
E\b_\sigma \times
E\b_{(\sigma,\partial \sigma)}
\longto
F\b_\sigma \times F\b_\sigma [-2]
\longto 
F\b_\sigma [-2] \longto
A\b[-2n]\;.
$$
Passing to the quotients modulo~$\mm_{A}$, we obtain~(5.3.1), with 
$\mu \: \quer F\b_\sigma [-2] \to \quer F\b_\sigma$ being induced 
by multiplication with the function $g \in \A\b(\Sigma)$.\qed

\bigskip
\noindent
{\bf 5.3b Lemma.} {\it If Poincar\'e duality holds for affine fans 
$\langle \sigma \rangle$
with a cone $\sigma$ of dimension $\le n$, then it also holds for every
quasi-convex fan $\Delta$ in dimension $n$.}

\medskip
\noindent
{\it Proof:}
To simplify notation, we introduce the abbreviation $\widetilde{A}\b := 
A\b[-2n]$. We embed the ``global''  duality homomorphism 
$$
\Phi \: E\b_\Delta \longto
\Hom_{A\b}(E\b_{(\Delta, \partial
\Delta)}, \widetilde{A}\b)
$$
induced by the pairing (5.1.3) into a commutative diagram of the 
following form:
$$
\klein
\def\longto{\kern-3pt\longrightarrow\kern-3pt}
\diagram
{0 & \longto & E\b_\Delta & \longto & C^0(\Delta, \partial \Delta; \E\b) 
& \longto & C^1(\Delta, \partial \Delta; \E\b) \cr  
\hbox{(5.3.3)}\kern-12pt & & \mapdown{\Phi} & & \mapdown{\Psi} & 
& \mapdown{\Theta} \cr 
0 & \longto & \Hom(E\b_{(\Delta, \partial \Delta)}, \widetilde A\b) &
\buildrel\kappa\over\longto & \bigoplus\limits_{\sigma \in \Delta^n}
   \Hom(E\b_{(\sigma,\partial\sigma)}, \widetilde{A}\b) & 
\buildrel\lambda\over\longto & \kern-5pt \bigoplus\limits_{\tau \in 
   {\Rond{\Delta}}^{_{n-1}}} \kern-5pt  \Hom(
   E\b_{(\tau,\partial\tau)},\widetilde A\b_\tau[2])
\cr}
$$
\smallskip
\noindent
where $\Hom$ abbreviates $\Hom_{A\b}$
and $\Psi, \Theta$ are the duality
homomorphisms corresponding to the dual
pairings $E\b_\sigma \times
E\b_{(\sigma, \partial \sigma)}
\to E\b_{(\sigma, \partial \sigma)}
\to \tilde A\b$ for $\sigma \in
\Delta^n$ and
$E\b_\tau \times
E\b_{(\tau, \partial \tau)}
\to E\b_{(\tau, \partial \tau)}
\to \tilde A\b_\tau[2]$ with suitably chosen evaluation maps. 
The upper row is exact,
while the lower one is a complex with an
injection $\kappa$. Since $\Psi$ and
$\Theta$ turn out to be isomorphisms, a
simple diagram chase yields that the same
holds for~$\Phi$ which proves the theorem.

The upper sequence is exact, since $\Delta$ is quasi-convex, see
Theorem~3.8. Now the
evaluation map $\epsilon \: 
E\b_{(\Delta,
\partial \Delta)} \to \widetilde A\b$
induces a system 
$(\epsilon_\sigma)_{\sigma \in \Delta^n}$
of evaluation maps  $\epsilon_\sigma\:
E\b_{(\sigma, \partial \sigma)} 
\subset E\b_{(\Delta, \partial \Delta)}
\to 
\widetilde A\b$:We have to show
$\epsilon_\sigma \not= 0$ for all
$\sigma \in \Delta$ resp. that the
homomorphism
$E\b_{(\sigma, \partial
\sigma)} \to E\b_{(\Delta, \partial
\Delta)}$ induces a non-trivial map
$\RR \cong \quer E^{2n}_{(\sigma, \partial
\sigma)} \to \quer E^{2n}_{(\Delta,
\partial \Delta)} \cong \RR$. In order
to do that we embed $\Delta$ into a
complete fan $\tilde \Delta$ and show
the corresponding fact for
$\quer E^{2n}_{(\sigma, \partial \sigma)}
\to \quer E^{2n}_{\tilde \Delta}$. But if
$\quer E^{2n}_{(\sigma, \partial \sigma)}
\to \quer E^{2n}_{\tilde \Delta}$ is
non-trivial, so is $\quer
E^{2n}_{(\sigma, \partial \sigma)} \to
\quer E^{2n}_{(\Delta, \partial
\Delta)}$, i.e., we may assume that
$\Delta$ is complete.  \par
To 
the quasi-convex fan $\Delta_0 :=
 \Delta \setminus \{\sigma\}$ 
corresponds an exact sequence 
$$
0 \wideto E\b_{(\Delta,
\Delta_0)} \cong E\b_{(\sigma, \partial 
\sigma)} \wideto E\b_{\Delta}
\wideto E\b_{\Delta_0} \wideto
0\;. 
$$ 
In the induced exact sequence
$\quer E^{2n}_{(\sigma, \partial 
\sigma)} \to \quer E^{2n}_{\Delta}
\to \quer  E^{2n}_{\Delta_0}$, 
the last term vanishes according to
4.4, since $\Delta_0$ is not 
complete. That yields our claim.   
\par
The system
of duality isomorphisms $E\b_\sigma
\to 
 \Hom(E\b_{(\sigma,\partial\sigma)}, 
\widetilde A\b)$ induced by the
compositions of the intersection
pairings with the $\epsilon_\sigma, \sigma
\in \Delta^n$, provides the isomorphism
$\Psi$. --- The map $\kappa$ associates
to a  homomorphism $\varphi \:
E\b_{(\Delta, \partial \Delta)} \to
\widetilde A\b$ its restrictions to the 
submodules $E\b_{(\sigma,\partial
\sigma)}$ of $E\b_{(\Delta, \partial
\Delta)}$. It is injective, since 
$\bigoplus_{\sigma  \in \Delta^n}
E\b_{(\sigma,\partial \sigma)} \cong
E\b_{(\Delta,  \Delta^{\le n-1})}$ is a
submodule of maximal rank in
$E\b_{(\Delta, \partial \Delta)}$:
For $h:=\prod_{\tau \in
\buildrel \circ \over \Delta^{n-1}}
h_\tau$, where $h_\tau \in A^2 \setminus
\{ 0 \}$ with kernel $V_\tau$, we have  $$
h E\b_{(\Delta, \partial\Delta)}
\subset E\b_{(\Delta, \Delta^{\le  n-1})}\ .
$$
This ends the discussion of the 
first rectangle in~(5.3.3).

The map $\lambda$ is composed of 
``restriction homomorphisms''
$$
\lambda_\tau^\sigma\:\Hom(E\b_{(\sigma,\partial \sigma)}, \widetilde A\b) 
\wideto 
\Hom(E\b_{(\tau,\partial \tau)}, \widetilde A\b_\tau[2]), \qquad \phi \mapsto 
\phi_\tau\;,
$$
where $\tau \prec_1 \sigma$ is a facet of $\sigma \in \Rond\Delta$.
In order to define them, we fix a
euclidean norm on $V$ and thus also  on
$V^* \cong A^2$. Now let $h_\tau \in A^2$
be the unique linear form of length $1$
that vanishes on $V_\tau$ and is
positive  on~${\Rond \sigma}$. Then we
use three exact sequences, where the 
first one is obvious: 
$$ 
0 \wideto
E\b_{(\sigma,\partial \sigma)} \wideto
E\b_\sigma \wideto  E\b_{\partial \sigma}
\wideto 0 \;.\leqno(5.3.4) 
$$ The
second one is multiplication with
$h_\tau$ and the projection onto its
cokernel: 
$$ 
0 \wideto \widetilde A\b~
{\buildrel \mu(h_\tau) \over \longto}~ 
 \widetilde A\b[2] \wideto \widetilde A\b_\tau[2] \wideto 0\;.
\leqno(5.3.5)
$$
Eventually the subfan $\partial_\tau \sigma :=  \partial\sigma 
\setminus \{\tau\}$ of $\partial\sigma$ gives the exact sequence
$$
0 \wideto E\b_{(\tau, \partial \tau)} \wideto E\b_{\partial \sigma}
\wideto E\b_{\partial_\tau  \sigma}
\wideto 0\;.\leqno(5.3.6)
$$
The associated $\Hom$-sequences provide 
a diagram

$$\klein
\def\longto{\kern-6pt\longrightarrow\kern-6pt}
\diagram{
    &   &   &   & \Ext(E_{\partial_\tau\sigma}\b,\widetilde A\b) 
                    &   &   \cr 
{\kern-50pt(5.3.7)}
    &   &   &   & \mapdown{}  
                    &   &   \cr 
  \Hom(E\b_\sigma, \widetilde A\b) & \longto & 
          \Hom(E\b_{(\sigma,\partial\sigma)}, \widetilde A\b) & 
              \buildrel \alpha\over\longto &
                  \Ext(E\b_{\partial\sigma}, \widetilde A\b)  
                    &   &   \cr 
    &   &   &   & \mapdown{\beta}  
                    &   &   \cr 
  \Hom( E_{(\tau,\partial\tau)}\b,  \widetilde A\b[2]) & \longto & 
          \Hom( E_{(\tau,\partial\tau)}\b,  \widetilde A_\tau\b[2]) & 
              \buildrel \gamma\over\longto &
                  \Ext(E_{(\tau,\partial\tau)}\b,  \widetilde A\b) 
                    & \longto & 
                          \Ext(E_{(\tau,\partial\tau)}\b, \widetilde A\b[2])  
\cr 
}
$$
with $\Ext=\Ext^1_{A\b}$. We show that $\gamma$ is an isomorphism; we 
then may set 
$\lambda_\tau^\sigma := 
\gamma^{-1}\circ\beta\circ\alpha$. Actually 
$\Ext^1_{A\b}(E\b_{(\tau, \partial \tau)},\widetilde A\b) \to 
\Ext^1_{A\b}(E\b_{(\tau, \partial \tau)}, \widetilde A\b[2])$ is the 
zero homomorphism, since it is induced by 
multiplication with $h_\tau$, which annihilates $E\b_{(\tau, \partial 
\tau)}$. On the other hand, the $A\b_\tau$-module 
$E\b_{(\tau, \partial \tau)}$ is a torsion 
module over $A\b$, so that 
$\Hom(E_{(\tau,\partial\tau)}\b, \widetilde A\b[2])$ vanishes.
\par

An explicit description of $\lambda^\sigma_\tau$ is as follows:
Let $\varphi \: E\b_{(\sigma, \partial \sigma)} \to 
\widetilde A\b$ be a homomorphism. Its ``restriction'' 
$\lambda^\sigma_\tau (\varphi)$ is the homomorphism
$\phi_\tau \: E\b_{(\tau, \partial \tau)} \to 
\widetilde A\b_\tau[2]$ given by this rule: For 
$g \in E\b_{(\tau, \partial \tau)}$, we choose a 
section $\hat g \in E\b_\sigma$ such that 
$\hat g|_{\partial \sigma}$ is the trivial extension of $g$ to 
$\partial \sigma$; then we have
$\phi_\tau (g) = \phi (h_\tau \hat g)|_\tau$.

Let us consider two different cones $\sigma = \sigma_1, \sigma_2 \in 
\Delta^n$ with intersection $\tau \in {\Rond \Delta}^{_{n-1}}$.
From the description of $\phi_\tau$, one easily derives that the 
compositions
$$
\Hom(E\b_{(\Delta, \partial \Delta)}, \widetilde A\b) \to \Hom 
(E\b_{(\sigma_{i}, 
\partial\sigma_{i})},\widetilde A\b) \to \Hom(E\b_{(\tau, \partial 
\tau)},\widetilde A\b_\tau[2])\leqno(5.3.8)
$$
coincide. --- Hence, we may define $\lambda$ in the usual way of a 
\v{C}ech coboundary homomorphism, starting from the appropriate 
$\lambda_\tau^\sigma$'s. Thus the lower row 
of  diagram (5.3.3) is a complex.

For the definition of $\Theta$, we need compatible evaluation 
homomorphisms 
$$
\epsilon_\tau\: E\b_{(\tau, \partial \tau)}\to \widetilde A\b_\tau[2]
$$ 
for $\tau \in {\Rond\Delta}^{_{n-1}}$. 
We choose them as the 
restrictions of the given evaluation map 
$E\b_{(\sigma, \partial \sigma)} \to \widetilde A\b$ 
for a $\sigma$ that includes $\tau$. As we have seen in (5.3.8), 
$\epsilon_\tau$ does not depend on the particular choice of $\sigma$. We 
still have to verify that $\epsilon_\tau$ is not the zero
homomorphism, i.e., we have to see that
$\lambda_\tau^\sigma$ is injective in
degree 0. In  diagram (5.3.7), we
have to show that $\alpha$ and $\beta$
are  injective in degree 0. By 4.4,
the vector spaces  $\quer E_\sigma^q$
vanish  for $q \ge 2n$; hence
$E\b_\sigma$ can be generated by elements
of degree $<2n$, and that yields the 
vanishing of $\Hom(E_\sigma\b; \widetilde
A\b)$ in degree 0.  According to
Lemma 5.4 below,  the exact sequence  
$$ 0
\wideto E\b_{(\sigma, \partial_\tau 
\sigma)} \wideto E\b_\sigma \wideto 
E\b_{\partial_\tau  \sigma} \wideto 0 
$$
is a free resolution of
$E\b_{\partial_\tau  \sigma}$, in
particular,  the module $\Ext^1
(E\b_{\partial_\tau  \sigma}, \widetilde
A\b)$ is a  quotient of
$\Hom(E\b_{(\sigma, \partial_\tau
\sigma)},\widetilde A\b)$,  and the
latter is trivial in degree $0$ by
5.4 below. --- For  $\tau \in
\Rond{\Delta}^{_{n-1}}$, the evaluation
homomorphisms  $\epsilon_\tau$ induce
isomorphisms $$ E\b_\tau \cong \Hom
(E\b_{(\tau,\partial \tau)}, \widetilde
A\b_\tau[2])\ , $$ which constitute the
isomorphism $\Theta$. Finally the
commutativity of the second square in
(5.3.3) follows from the above
explicit description of the restriction
homomorphisms $\lambda^\sigma_\tau$ and
the appropriate choice of the evaluation
homomorphisms $\varepsilon_\tau$. This
ends the discussion of (5.3b) and
the proof of the theorem. \qed

\medskip
For the notation used in the following
result that has been used in the proof of
5.3b,  we refer to~(0.D.2). 
\medskip

\noindent
{\bf 5.4  Proposition.} {\sl Let $\sigma$ be a cone of 
dimension $n$ and $\Lambda \subset \partial \sigma$ be a fan such that 
$\pi ( \Lambda)$ is a quasi-convex subfan of $\Lambda_\sigma$. Then 
$E\b_{(\sigma, \Lambda)}$ is a free 
$A\b$-module, and, if in
addition $\Lambda $ is a  proper subfan,
$\quer E^q_{(\sigma, \Lambda)}=0$ for $q
\ge 2n$.} 

\medskip \noindent 
{\sl Proof.} As in 0.D, we choose a line $L \subset V$ intersecting $\Rond
\sigma$ and set $B\b :=\pi^*( S\b((V/L)^*))) \subset A\b$; 
furthermore, we write
$A\b = B\b[T]$ with a linear form 
$T \in A^2$. The exact sequence of 
$A\b$-modules
$$
0 \wideto E\b_{(\sigma, \Lambda)} \wideto E\b_\sigma \wideto E\b_{\Lambda} 
\wideto 0\
$$
induces an exact $\Tor$-sequence
$$
\Tor_2^{A\b}(E\b_\Lambda, \RR\b) \wideto \Tor_1^{A\b}(E\b_{(\sigma, 
\Lambda)},\RR\b) \to 0 \to \Tor_1^{A\b} (E\b_\Lambda, \RR\b) \wideto 
\quer E\b_{(\sigma, \Lambda)} \to \quer E\b_\sigma
$$
since $E\b_\sigma$ is free. If $\Tor_2^{A\b}(E\b_\Lambda, \RR\b)=0$, then 
also $\Tor_1^{A\b}(E\b_{(\sigma, \Lambda)},\RR\b)$, and $E\b_{(\sigma, 
\Lambda)}$ is a free $A\b$-module by 
section 0.B. Since the fan
$\langle \sigma \rangle$ is not  complete,
$\quer E^q_{\sigma}=0$ for $q \ge 2n$; if
the same holds for  $\Tor_1^{A\b}
(E\b_\Lambda, \RR\b)$, then it follows
for $\quer  E^q_{(\sigma, \Lambda)}$ as
well. It thus remains to determine 
$\Tor^{A\b}_i(E\b_{\Lambda},\RR\b)$.
Once again we use the exact sequence 
$$ 
0 \wideto
\RR\b[T][-2] ~{\buildrel \mu(T) \over
\longto}~ \RR\b[T]  \wideto \RR\b \wideto
0 
$$ 
of $A\b$-module homomorphisms of
degree 0, where  $\RR\b[T]$ denotes the
$A\b$-module $A\b/(\mm_{B\b} A\b)
=B\b / \mm_{B\b}[T]$, and $
\mm_{B\b} :=B^{>0}$, the maximal
homogeneous ideal of $B\b$. Using the
identity  
$$
  \Tor_i^{B\b}(E\b_\Lambda, \RR\b) \widecong 
  \Tor_i^{B\b[T]}(E\b_\Lambda, \RR\b[T]) 
$$
we obtain the exact sequence
$$
\dots \to \Tor_1^{B\b}(E\b_\Lambda, \RR\b[-2]) 
\wideto \Tor_1^{B\b}(E\b_\Lambda, \RR\b) \to 
\Tor_1^{A\b}(E\b_\Lambda, \RR\b) 
$$
$$
\to 
E\b_\Lambda\otimes_{B\b}\RR\b [-2] \to 
E\b_\Lambda\otimes_{B\b}\RR\b \;,
$$
where $\Tor_i^{B\b}(E\b_\Lambda, \RR\b)$ 
vanishes for $i \ge 1$ since 
$E\b_\Lambda$ is a free $B\b$-module. 
This yields the desired description:
$$
\Tor_i^{A\b}(E\b_\Lambda, \RR\b)~ =~ 
\cases{\ker \Bigl(\mu(T) \:  
E\b_\Lambda\otimes_{B\b}\RR\b [-2] 
\wideto  
E\b_\Lambda\otimes_{B\b}\RR\b \Bigr), & if
$i=1$; \cr  0, & if $i \ge 2\;.$ \cr}   
$$
Eventually, if $\pi (\Lambda) \subset \Lambda_\sigma$ is not complete, 
then the vector space 
$E\b_\Lambda\otimes_{B\b}\RR\b [-2]$
vanishes in degrees $\ge 2n$;  hence, the
same  holds for $\Tor_1^{A\b}
(E\b_\Lambda, \RR\b)$.\qed

\medskip
\noindent 
{\bf 5.5 Remark.} For every purely $n$-dimensional 
fan $\Delta$ we can define an evaluation
map $E\b_{(\Delta, \partial \Delta)}
\to A\b[-2n]$ as the restriction of an
evaluation map $E\b_{\tilde \Delta}
\to A\b[-2n]$ for a ``completion''
$\tilde \Delta$ of $\Delta$. It
provides a homomorphism $E\b_\Delta \to 
\Hom(E\b_{(\Delta, \partial\Delta)}, A\b[-2n])$ 
via the intersection pairing. In
accordance with  the proof of
5.3b, that is an isomorphism
whenever $\tilde H^0(\Delta, 
\partial\Delta; \E\b) = 0$, or
equivalently (see~3.7), if  $\tilde
H^0(\Delta_\sigma, \partial
\Delta_\sigma; \RR\b) = \{0\}$ for each cone 
$\sigma \in \Delta$. In more geometrical
terms,  $\Delta$ has to be both facet-connected and
locally facet-connected,  where we call a fan {\it
locally facet-connected\/} if for  each non-zero
cone $\sigma \in \Delta$, the fan 
$\Delta_\sigma$ is facet-connected.

The smallest example of a three-dimensional fan that is both 
facet-connected and locally facet-connected, but not quasi-convex, 
is provided by the fan swept out by the ``vertical'' facets of a 
triangular prism.
\medskip
Since the dual pairing $E\b_\Delta \times E\b_{(\Delta,
\partial \Delta)} \to
A\b[-2n]$ of $A\b$-modules induces a dual pairing of real vector 
spaces $\quer E\b_\Delta \times \quer E\b_{(\Delta,
\partial \Delta)} \to
\RR\b[-2n]$, we obtain the following consequence.
\bigskip

\noindent
{\bf 5.6 Corollary.} {\it If $\Delta$ is a quasi-convex fan of 
dimension~$n$, then we have 
$$
\dim \quer E^q_\Delta \;=\; \dim \quer E^{2n-q}_{(\Delta,
\partial \Delta)}\; . \eqno\qed  
$$}

We finally are prepared to prove the ``Combinatorial Hard Lefschetz'' 
Theorem~4.6.

\bigskip
\noindent
{\it Proof of Theorem 4.6:} Since $f$ is
strictly convex, its graph $\Gamma_f$ in
$V \times \RR$ is the support of the
boundary fan $\partial \gamma$ of the
$(n+1)$-dimensional  cone $\gamma :=
\gamma^+(f)$ in $V \times \RR$
as we have  seen in 0.D. Denote  with
$\H\b$ a minimal extension sheaf on
$\partial \gamma$.  Then 
$\tau \mapsto
\H\b(f(\tau))$ is a minimal extension
sheaf on $\Delta$, that we identify with
$\E\b$. Then in analogy to (4.3.2) the
residue class module of the
$A\b[T]$-module $H\b_{\partial \gamma}$
satisfies
$$
\quer H\b_{\partial \gamma} \cong
\quer E\b_\Delta / f \quer E\b_\Delta
=
\coker (\quer\mu: \quer E\b_\Delta [-2]
\longto \quer E\b_\Delta)
$$
where $\quer E\b_\Delta = (A\b/\mm)
\otimes_{A\b} E\b_\Delta$. 
Now the vanishing  condition
$V(\gamma)$ yields the surjectivity of
$\quer \mu^{2q}$  for $2q \ge n-1$.
On the other hand the map $\mu$
is  selfadjoint with respect to the dual
pairing  $E\b_\Delta \times E\b_{\Delta}
\to A\b[-2n]$ as well as 
$\quer \mu$ with respect to  $\quer
E\b_\Delta \times \quer E\b_{\Delta} \to
\RR\b[-2n]$. Hence by Poincar\'e duality
the surjectivity of $\quer \mu^{2q}$ for
$2q \ge n-1$ implies the injectivity of
$\quer \mu^{2q}$ for $2q \le n-1$. \qed

\goodbreak\goodbreak\goodbreak\goodbreak\ 
\bigskip\goodbreak\goodbreak\goodbreak\goodbreak\ \medskip
\def\litemitem{\par\noindent
                   \hangindent=45pt\ltextindent}
\def\ltextindent#1{\hbox to \hangindent{#1\hss}\ignorespaces}
{\klein\font\xisc=cmcsc10 at 11 truept\def\sc{\xisc}
\vbox
{\centerline{\bf References}
\nobreak\nobreak\nobreak\nobreak\smallskip\noindent

\litemitem{[BBD]} {\sc A.~Beilinson, J.~Bernstein, P.~Deligne:} {\it 
Faisceaux pervers}, Ast\'erisque, vol.~{\bf 100}, Soc. Math. France, 
1982.}

\litemitem{[BBFK]} {\sc G.~Barthel, J.-P.~Brasselet, F.-H.~Fieseler, 
L.~Kaup:} {\it Equivariant Intersection Cohomology of Toric Varieties}, 
in: {\it Algebraic Geometry: Hirzebruch 70}. Contemp.\ Math.\ 
vol.~{\bf 241}, Amer. Math. Soc., Providence, R.I., 1999, 45--68.

\litemitem{[BeLu]} {\sc J.~Bernstein, V.~Lunts:} {\it Equivariant 
Sheaves and Functors}, Lecture Notes in Math. vol.~{\bf 1578}, 
Springer-Verlag Berlin etc., 1993.  

\litemitem{[BrMPh]} {\sc T.~Braden, R.~McPherson:} {\it 
Intersection homology of toric varieties and a conjecture of Kalai}, 
Comment. Math. Helv.~{\bf 74} (1999), 442--455. 

\litemitem{[BreLu$_1$]} {\sc P.~Bressler, V.~Lunts:} {\it Toric Varieties 
and Minimal Complexes}, (pr)e-print {\tt alg-geom/9712007}, 1997. 

\litemitem{[BreLu$_2$]} {\sc P.~Bressler, V.~Lunts:} {\it Intersection 
cohomology on non-rational polytopes}, (pr)e-print {\tt math.AG/0002006}, 
2000. 

\litemitem{[Bri]} {\sc M.~Brion:} {\it The Structure of the Polytope 
Algebra}, T\^ohoku Math.\ J.~{\bf 49} (1997), 1--32.

\litemitem{[Fi$_1$]} {\sc K.-H.~Fieseler:} {\it Rational Intersection 
Cohomology of Projective Toric Varieties}, J.~reine angew.\ Math.\ 
(Crelle)~{\bf 413} (1991), 88--98. 

\litemitem{[Fi$_2$]} {\sc K.-H.~Fieseler:} {\it Towards a Combinatorial 
Intersection Cohomology for Fans}, to appear in Comptes Rendues Acad.\ 
Sc.\ Paris, 2000. 

\litemitem{[GoKoMPh]} \quad{\sc M.~Goresky, R.~Kottwitz, R.~MacPherson:} 
{\it Equivariant Cohomology, Koszul Duality, and the Localization Theorem}, 
Invent. Math. {\bf 131} (1998), 25--83. 

\litemitem{[Ish]} {\sc M.-N.~Ishida:} {\it Torus Embeddings and Algebraic 
Intersection Complexes~I, II\/}, (pr)e-print {\tt alg-geom/9403008, 
alg-geom/9403009}, 1993. Revised Version: {\it Combinatorial and Algebraic 
Intersection Complexes of Toric Varieties\/}, available from the
author's web site, {\tt http://www.math.tohoku.ac.jp/~ishida/}.

\litemitem{[MCo]} {\sc M.~McConnell:} {\it Intersection Cohomology 
of Toric Varieties}, preprint (available on the author's web site, 
{\tt http://www.math.okstate.edu/\~{}mmcconn/}), 1997. 

\litemitem{[Oda]} {\sc T.~Oda:} {\it The Intersection Cohomology and 
Toric Varieties}, in: {\sc T.~Hibi} (ed.): {\it Modern Aspects of 
Combinatorial Structure on Convex Polytopes}, RIMS {\it Kokyu\-roku} 
{\bf 857} (Jan. 1994), 99--112.

\litemitem{[S]} {\sc R.~Stanley:} {\it Generalized h-vectors,
intersection cohomology of toric varieties and related results
}, in {\sc M.~Nagata, H.~Matsumura} (eds.): {\it Commutative Algebra 
and Combinatorics},  Adv.\ Stud.\ Pure Math.~{\bf 11}, Kinokunia, Tokyo, 
and North Holland, Amsterdam/New York, 1987, 187-213.
\par
}

\bigskip\bigskip\goodbreak\noindent

\vbox{
{\klein Addresses of authors}
\bigskip
\nobreak
\hbox to \hsize{\hfill
{\klein
\parskip 0pt
\baselineskip=12pt
\vtop {\vskip 0.5truecm} {\hskip 0.1truecm}
\vbox{\hsize=0.5\hsize
\obeylines{G. Barthel, L. Kaup
Fachbereich Mathematik
\quad und Statistik
Universit\"at Konstanz
Fach D 203
D-78457 Konstanz
e-mail: 
{\tt Gottfried.Barthel@uni-konstanz.de
Ludger.Kaup@uni-konstanz.de}
}}\hfill
\vbox{\hsize=0.45\hsize
\obeylines{K.-H. Fieseler
Matematiska Institutionen
Box 480
Uppsala Universitet
SE-75106 Uppsala
e-mail: {\tt khf@math.uu.se}
{\ }
J.P. Brasselet.
IML/CNRS, Luminy Case 907
F-13288 Marseille Cedex 9
e-mail: {\tt jpb@iml.univ-mrs.fr}}}\hfill}}
}
\bye